\input amstex
\input eplain
\frenchspacing
\documentstyle{amsppt}
\magnification=\magstep1
\vsize=19.5cm
\NoBlackBoxes
\def\mmod{\mathop{\text{\rm mod}}\nolimits^*}
\def\Tr{{\text{\rm Tr}}}
\def\ab{{\text{\rm ab}}}
\def\im{{\text{\rm im}}}
\def\Fr{\mathop{\text{\rm Frob}}\nolimits}
\def\Art{\psi}
\def\Hil{\mathop{\text{\rm Hil}}\nolimits}
\def\Gal{{\text{\rm Gal}}}
\def\GL{{\text{\rm GL}}}
\def\PGL{{\text{\rm PGL}}}
\def\Sl{{\text{\rm SL}}}
\def\Cl{{\text{\rm Cl}}}
\def\Hom{{\text{\rm Hom}}}
\def\Aut{{\text{\rm Aut}}}
\def\End{{\text{\rm End}}}
\def\tor{{\text{\rm tor}}}
\def\ker{\mathop{\text{\rm ker}}}
\def\ord{{\text{\rm ord}}}
\def\sign{{\text{\rm sign}}}

\def\tto{\longrightarrow}

\def\mapright^#1{\ \smash{\mathop{\longrightarrow}\limits^{#1}}\ }
\def\isar{\ \smash{\mathop{\longrightarrow}\limits^\sim}\ }
\def\mod{\bmod}
\def\congr{ \equiv }
\def\iso{ \cong }

\def\AA{\bold A}
\def\Que{\bold Q}
\def\CC{\bold C}
\def\FF{\bold F}
\def\HH{\bold H}
\def\PP{\bold P}
\def\RR{\bold R}
\def\Zee{\bold Z}
\def\gotha{{\goth a}}
\def\gothb{{\goth b}}
\def\gothc{{\goth c}}
\def\gothf{{\goth f}}
\def\gothm{{\goth m}}
\def\gothn{{\goth n}}
\def\gothp{{\goth p}}
\def\gothq{{\goth q}}

\def\endexample{\hfill$\lozenge$}
\newcount\refCount
\def\newref#1 {\advance\refCount by 1
\expandafter\edef\csname#1\endcsname{\the\refCount}}
\newref AT   
\newref BG   
\newref BS   
\newref BW   
\newref CF   
\newref COA  
\newref COB  
\newref CX   
\newref DE   
\newref FI   
\newref HS   
\newref KI   
\newref LA   
\newref LAEF 
\newref PO   
\newref SA   
\newref LS   
\newref SC   
\newref SE   
\newref Shi  
\newref SI   
\newref SIB  
\newref ST   
\newref STB  
\newref SL   
\newref VO   
\newref WA   
\newref WE   
\topmatter
\title 
Computational class field theory
\endtitle
\author Henri Cohen and Peter Stevenhagen\endauthor
\address
Laboratoire A2X, U.M.R. 5465 du C.N.R.S.,
Universit\'e Bordeaux I, 351 Cours de la Lib\'eration,
33405 Talence Cedex, France
\endaddress
\email cohen\@math.u-bordeaux1.fr\endemail

\address Mathematisch Instituut,
Universiteit Leiden, Postbus 9512, 2300 RA Leiden, The Netherlands
\endaddress
\email psh\@math.leidenuniv.nl\endemail
\abstract
Class field theory furnishes an intrinsic description of the abelian
extensions of a number field that is in many cases not of an
immediate algorithmic nature.
We outline the algorithms available for the explicit
computation of such extensions.
\endabstract
\endtopmatter

\document

\head 1. Introduction
\endhead

\noindent
Class field theory is a 20th century theory
describing the set of finite {\it abelian\/} extensions $L$ of
certain base fields $K$ of arithmetic type.
It provides a canonical description of the Galois groups $\Gal(L/K)$
in terms of objects defined `inside $K$', and gives rise to an explicit
determination of the maximal abelian quotient $G_K^{\ab}$ of the
absolute Galois group $G_K$ of $K$.
In the classical examples $K$ is either a {\it global field}, i.e., a number
field or a function field in one variable over a finite field,
or a {\it local field\/} obtained by completing a global field at one
of its primes.
In this paper, which takes an algorithmic approach,
we restrict to the fundamental case in which
the base field $K$ is a number field.
By doing so,
we avoid the complications arising for $p$-extensions
in characteristic $p>0$.

The description of $G_K^{\ab}$ for a number field $K$ that is
provided by class field theory
can be seen as a first step towards a complete description of the
full group $G_K\subset G_\Que$.
At the moment, such a description is still far away, and it is not even
clear what kind of description one might hope to achieve.
Grothendieck's anabelian Galois theory and his theory of
{\it dessins d'enfant\/} [\LS] constitute one direction of progress, and the
largely conjectural {\it Langlands program\/} [\BG] provides an other
approach.
Despite all efforts and partial results [\VO],
a concrete question such as the {\it inverse problem
of Galois theory\/}, which asks whether for a number field $K$,
all finite groups $G$ occur as the Galois group of some finite extension $L/K$,
remains unanswered for all $K$.


A standard method to gain insight into the structure of $G_K$, and
to realize certain types of Galois groups over $K$ as quotients of $G_K$,
consists in studying the action of $G_K$ on `arithmetical objects'
related to $K$, such as the division points in $\overline \Que$ of
various algebraic groups defined over $K$.
A good example is the {\it Galois representation\/}
arising from the group $E[m](\overline \Que)$ of $m$-torsion points 
of an elliptic curve $E$ that is defined over~$K$.
The action of $G_K$ on $E[m](\overline \Que)$ factors via a finite quotient
$T_m\subset \GL_2(\Zee/m\Zee)$ of $G_K$, and
much is known [\SE] about the groups~$T_m$.
Elliptic curves with {\it complex multiplication\/} by an order
in an imaginary quadratic field $K$ give rise to {\it abelian\/}
extensions of $K$, and yield a particularly explicit
instance of class field theory.

For the much simpler example of the multiplicative group ${\bold G}_m$,
the division points of ${\bold G}_m(\overline \Que)$
are the {\it roots of unity\/} in $\overline \Que$.
The extensions of $K$ they generate are the {\it cyclotomic extensions\/}
of $K$.
As the Galois group of the extension $K\subset K(\zeta_m)$ obtained
by adjoining a primitive $m$-th root of unity $\zeta_m$ to $K$
naturally embeds into $(\Zee/m\Zee)^*$, all cyclotomic extensions
are abelian.
For $K=\Que$, Kronecker discovered in 1853 that {\it all\/}
abelian extensions are accounted for in this way.
\proclaim
{1.1. Kronecker-Weber theorem}
Every finite abelian extension $\Que\subset L$ is contained in some
cyclotomic extension $\Que\subset \Que(\zeta_m)$.
\endproclaim\noindent
Over number fields $K\ne\Que$, there are more abelian extensions
than just cyclotomic ones, and the analogue of 1.1 is what 
class field theory provides: every abelian extension $K\subset L$
is contained in some {\it ray class field extension\/} $K\subset H_\gothm$.
Unfortunately, the theory does not provide
a `natural' system of generators for the fields $H_\gothm$ to play
the role of the roots of unity in 1.1.
Finding such a system for all $K$ is one of the
Hilbert problems from 1900 that is still open.
Notwithstanding this problem,
class field theory is in principle constructive, and once one finds
in some way a possible generator of $H_\gothm$ over $K$,
it is not difficult to verify that it does generate $H_\gothm$.
The information on $H_\gothm$ we have is essentially an
intrinsic description, in terms of the splitting and ramification 
of the primes in the extension $K\subset H_\gothm$, of
the Galois group $\Gal(H_\gothm/K)$
as a {\it ray class group\/} $\Cl_\gothm$.
This group replaces the group $(\Zee/m\Zee)^*$ that occurs implicitly in 1.1
as the underlying Galois group:
$$
\eqalign{
(\Zee/m\Zee)^*&\isar \Gal(\Que(\zeta_m)/\Que)\cr
(a\mod m)&\tto (\sigma_a: \zeta_m\mapsto \zeta_m^a).\cr
}
\leqno{(1.2)}
$$
We can in principle find generators for any specific class
field by combining our knowledge of its ramification data with a
classical method to generate arbitrary solvable field extensions:
the adjunction of {\it radicals}.
More formally, we call an extension $L$ of an arbitrary field $K$ a
{\it radical extension\/} if $L$ is contained in the splitting field over $K$
of a finite collection of polynomials of the form $X^n-a$,
with $n\in\Zee_{\ge1}$ not divisible by $\text{char}(K)$ and $a\in K$.
If the collection of polynomials can be chosen such that
$K$ contains a primitive $n$-th root of unity for each polynomial $X^n-a$
in the collection, then the radical extension $K\subset L$ is said to be
a {\it Kummer extension\/}.
Galois theory tells us that every Kummer extension
is abelian and, conversely, that an abelian extension $K\subset L$
of exponent $n$ is Kummer in case $K$ contains a primitive $n$-th root of unity.
Here the {\it exponent\/} of an abelian extension $K\subset L$
is the smallest positive integer $n$ that annihilates $\Gal(L/K)$.
Thus, for every finite abelian extension $K\subset L$ of a number field $K$,
there exists a cyclotomic extension $K\subset K(\zeta)$ such that
the `base-changed' extension $K(\zeta)\subset L(\zeta)$ is Kummer.

In Section 5, we compute the class fields of $K$ as subfields of
Kummer extensions of $K(\zeta)$ for suitable cyclotomic extensions
$K(\zeta)$ of $K$.
The practical problem of the method is that the auxiliary fields $K(\zeta)$
may be much larger than the base field $K$,
and this limits its use to not-too-large examples.

If $K$ is imaginary quadratic, elliptic curves with complex multiplication
solve the Hilbert problem for $K$, and this yields methods that
are much faster than the Kummer extension constructions for general $K$.
We describe these complex multiplication methods in some detail
in our Sections~6 to 8.
We do not discuss their extension to abelian varieties
with complex multiplication [\Shi], or the analytic generation of class fields
of totally real number fields $K$ using {\it Stark units\/}
[\COB, Chapter~6].

\head 2. Class field theory
\endhead

\noindent
Class field theory generalizes 1.1 by focusing on
the Galois group $(\Zee/m\Zee)^*$ of the
cyclotomic extension $\Que\subset \Que(\zeta_m)$
rather than on the specific generator $\zeta_m$.
The extension $\Que\subset \Que(\zeta_m)$ is unramified
at all primes $p\nmid m$, and the splitting behavior
of such $p$ only depends on the residue class $(p\mod m)\in (\Zee/m\Zee)^*$.
More precisely, the residue class degree $f_p=[\FF_p(\zeta_m):\FF_p]$ of
the primes over $p\nmid m$ equals the order of the {\it Frobenius automorphism\/}
$(\sigma_p: \zeta_m\mapsto\zeta_m^p)\in \Gal(\Que(\zeta_m)/\Que)$,
and this is the order of $(p\mod m)\in (\Zee/m\Zee)^*$ under the
standard identification 1.2.

Now let $K\subset L$ be {\it any\/} abelian extension of number fields.
Then we have [\ST, Section 15],
for each prime $\gothp$ of $K$ that is unramified in $L$, a
unique element $\Fr_\gothp\in\Gal(L/K)$ that induces the Frobenius
automorphism $x \mapsto x^{\#k_\gothp}$ on the residue class field
extensions $k_\gothp\subset k_\gothq$ for the primes $\gothq$ in $L$
extending $\gothp$.
The order of this {\it Frobenius automorphism\/} $\Fr_\gothp$ of $\gothp$
in $\Gal(L/K)$ equals the residue class degree $[k_\gothq:k_\gothp]$,
and the subgroup $\langle \Fr_\gothp\rangle \subset \Gal(L/K)$ is
the decomposition group of $\gothp$.

We define the {\it Artin map\/} for $L/K$ as the homomorphism
$$
\eqalign{
\psi_{L/K}:\quad I_K(\Delta_{L/K})&\longmapsto \Gal(L/K)\cr
\gothp       &\longmapsto \Fr_\gothp\cr}
\leqno{(2.1)}
$$
on the group $I_K(\Delta_{L/K})$ of fractional $\Zee_K$-ideals generated by 
the primes $\gothp$ of $K$ that do not divide the discriminant
$\Delta_{L/K} $ of the extension $K\subset L$. 
Such primes $\gothp$ are known to be unramified in $L$ by [\ST, Theorem 8.5].
For an ideal $\gotha\in I_K(\Delta_{L/K})$,
we call $\psi_{L/K}(\gotha)$ the {\it Artin symbol\/} of
$\gotha$ in $\Gal(L/K)$.

For $K=\Que$, we can rephrase 1.1 as follows.
\proclaim
{2.2. Kronecker-Weber theorem}
If $\Que\subset L$ is an abelian extension, there exists an integer
$m\in\Zee_{>0}$ such that the kernel of the Artin map 
$\psi_{L/\Que}$ contains all $\Zee$-ideals $x\Zee$ with $x>0$
and $x\congr 1\mod m$.
\endproclaim\noindent
The equivalence of 1.1 and 2.2 follows from the analytic fact that
an extension of number fields is trivial if all primes
outside a density zero subset split completely in it.
Thus, if all primes $p\congr1\mod m$ split completely
in $\Que\subset L$, then all primes of degree one are split in
$\Que(\zeta_m)\subset L(\zeta_m)$ 
and $L$ is contained in the cyclotomic field $\Que(\zeta_m)$.

The positivity condition on $x$ in 2.2 can be omitted if the
primes $p\congr-1\mod m$ also split completely in $L$, i.e., if
$L$ is totally real, and contained in the maximal real subfield
$\Que(\zeta_m+\zeta_m^{-1})$ of $\Que(\zeta_m)$.
The values of $m$ that one can take in 2.2 are the multiples of 
some minimal positive integer, the {\it conductor\/} of $\Que\subset L$.
It is the smallest integer $m$ for which $\Que(\zeta_m)$ contains $L$.
The prime divisors of the conductor are exactly the primes that
ramify in $L$, and $p^2$ divides the conductor
if and only if $p$ is {\it wildly\/} ramified in $L$.

For a quadratic field $L$ of discriminant $d$, the conductor equals $|d|$,
and 2.2 states that the Legendre symbol ${d\overwithdelims() x}$
only depends on $x$ modulo $|d|$.
This is Euler's version of the quadratic reciprocity law.
The main statement of class field theory is the analogue of 2.2
over arbitrary number fields $K$.
\proclaim{2.3. Artin's reciprocity law}
If $K\subset L$ is an abelian extension, there exists a non-zero
ideal $\gothm_0\subset\Zee_K$ such that the kernel of the Artin map
$\psi_{L/K}$ in $(2.1)$ contains all principal $\Zee_K$-ideals $x\Zee_K$ with
$x$ totally positive and $x\congr 1\mod \gothm_0$.
\endproclaim\noindent
This somewhat innocuous looking statement is a highly
non-trivial fact. 
It shows that there is a powerful global connection relating the splitting
behavior in $L$ of {\it different\/} primes of $K$.
Just like 2.2 implies the quadratic reciprocity law, 
Artin's reciprocity law implies the general
{\it power reciprocity laws\/} from from algebraic number theory
(see [\AT, Chapter 12, \S4] or [\CF, p. 353]).

It is customary to treat the positivity conditions at the real primes of $K$
and the congruence modulo $\gothm_0$ in 2.3 on equal footing.
To this end, one formally defines a {\it modulus\/} $\gothm$ of $K$ to 
be a non-zero $\Zee_K$-ideal $\gothm_0$ times a subset $\gothm_\infty$
of the real primes of $K$.
For a modulus $\gothm=\gothm_0\gothm_\infty$, we write 
$$
x\congr 1\mmod \gothm
$$
if $x$ satisfies $\ord_\gothp(x-1)\ge\ord_\gothp(\gothm_0)$
at the primes $\gothp$
dividing the {\it finite part\/} $\gothm_0$, and
$x$ is positive at the real primes in the {\it infinite part\/} $\gothm_\infty$
of $\gothm$.

In the language of moduli, 2.3 asserts that there exists a modulus $\gothm$
such that the kernel $\ker\psi_{L/K}$ of the Artin map
contains the {\it ray group\/} $R_\gothm$ of
principal $\Zee_K$-ideals $x\Zee_K$ generated by elements
$x\congr 1\mmod \gothm$.
As in the case of 2.2, the set of these {\it admissible\/} moduli for
$K\subset L$ consists of the multiples $\gothm$ of some minimal 
modulus $\gothf_{L/K}$, the {\it conductor\/} of $K\subset L$.
The primes occurring in $\gothf_{L/K}$ are the primes of $K$,
both finite and infinite, that ramify in $L$.
An infinite prime of $K$ is said to ramify in $L$ if it is real
but has complex extensions to $L$.
As for $K=\Que$, a finite prime $\gothp$ occurs with
higher multiplicity in the conductor if and only if it is
wildly ramified in $L$.

If $\gothm=\gothm_0\gothm_\infty$ is an admissible modulus for $K\subset L$, and 
$I_\gothm$ denotes the group of fractional $\Zee_K$-ideals generated
by the primes $\gothp$ coprime to $\gothm_0$, the Artin map
induces a homomorphism
$$
\eqalign{
\psi_{L/K}:\quad \Cl_\gothm=I_\gothm/R_\gothm  &\tto \Gal(L/K)\cr
[\gothp]      &\longmapsto \Fr_\gothp\cr}
\leqno{(2.4)}
$$
on the {\it ray class group\/} $\Cl_\gothm=I_\gothm/R_\gothm$ modulo $\gothm$.
Our earlier remark on the triviality of extensions in which
almost all primes split completely implies that it is {\it surjective\/}.
By the Chebotarev density theorem [\SL], even more is true:
the Frobenius automorphisms $\Fr_\gothp$ for $\gothp\in I_\gothm$
are {\it equidistributed\/} over the Galois group $\Gal(L/K)$.
In particular, a modulus $\gothm$ is admissible for an abelian
extension $K\subset L$ if and only if (almost) all primes $\gothp\in
R_\gothm$ of $K$ split completely in $L$.

As the order of the Frobenius automorphism
$\Fr_\gothp\in\Gal(L/K)$ equals the residue class
degree $f_\gothp$ of the primes $\gothq$ in $L$ lying over $\gothp$,
the norm $N_{L/K}(\gothq)=\gothp^{f_\gothp}$ of every prime ideal
$\gothq$ in $\Zee_L$ coprime to $\gothm$ is contained in the
kernel of the Artin map.
A non-trivial index calculation shows that the norms of the
$\Zee_L$-ideals coprime to $\gothm$ actually generate the
kernel in (2.4).
In other words, the {\it ideal group\/} $A_\gothm\subset I_\gothm$
that corresponds to $L$ in the sense that we have 
$\ker\psi_{L/K}=A_\gothm/R_\gothm$, is equal to
$$
A_\gothm=N_{L/K}(I_{\gothm\Zee_L})\cdot R_\gothm.
\leqno(2.5)
$$
The {\it existence theorem\/} from class field theory states that
for every modulus $\gothm$ of~$K$, there exists an extension
$K\subset L=H_\gothm$ for which the map $\psi_{L/K}$ in (2.4)
is an isomorphism.
Inside some fixed algebraic closure $\overline K$ of $K$, the extension 
$H_\gothm$ is uniquely determined as the maximal abelian extension $L$
of $K$ in which all primes in the ray group $R_\gothm$ split completely.
It is the {\it ray class field\/} $H_\gothm$ modulo $\gothm$ mentioned
in the introduction, for which the analogue of 1.1 holds over $K$. 
If $K\subset L$ is abelian, we have $L\subset H_\gothm$
whenever $\gothm$ is an admissible modulus for $L$.
For $L=H_\gothm$, we have $A_\gothm=R_\gothm$ in (2.5), and 
an Artin isomorphism $\Cl_\gothm\isar \Gal(H_\gothm/K)$.
\medskip\noindent
{\bf 2.6. Examples.}
1. 
It will not come as a surprise that
for $K=\Que$, the ray class field modulo $(m)\cdot\infty$ is the cyclotomic
field $\Que(\zeta_m)$, and the ray class group $\Cl_{(m)\cdot\infty}$
the familiar group $(\Zee/m\Zee)^*$ acting on the $m$-th roots of unity.
Leaving out the real prime $\infty$ of~$\Que$, we find
the ray class field modulo $(m)$ to be
the maximal real subfield $\Que(\zeta_m+\zeta_m^{-1})$ of $\Que(\zeta_m)$.
This is the maximal subfield in which the real prime $\infty$ is unramified.

2. The ray class field of conductor $\gothm= (1)$
is the {\it Hilbert class field\/} $H=H_1$ of $K$.
It is the largest abelian extension of $K$ that is unramified at
all primes of $K$, both finite and infinite.
As $I_1$ and $R_1$ are the groups of all fractional and all
principal fractional $\Zee_K$-ideals, respectively,
the Galois group $\Gal(H/K)$ is isomorphic to the ordinary class group
$\Cl_K$ of $K$,
and the primes of $K$ that split completely in $H$ are
precisely the {\it principal\/} prime ideals of $K$.
This peculiar fact makes it possible to derive information about the class
group of $K$ from the existence of unramified extensions of $K$, and conversely.
\endexample
\medskip\noindent
The ray group $R_\gothm$ is contained in the subgroup $P_\gothm\subset 
I_\gothm$ of principal ideals in $I_\gothm$,
and the quotient $I_\gothm/P_\gothm$
is the class group $\Cl_K$ of $K$ for all $\gothm$.
Thus, the ray class group $\Cl_\gothm=I_\gothm/R_\gothm$ is an extension
of $\Cl_K$ by a finite abelian group $P_\gothm/R_\gothm$
that generalizes the groups $(\Zee/m\Zee)^*$ from 1.2.
More precisely, we have a natural exact sequence
$$
\Zee_K^*\tto (\Zee_K/\gothm)^* \tto \Cl_\gothm \tto \Cl_K\tto 0
\leqno(2.7)
$$
in which the residue class of $x\in \Zee_K$ coprime to $\gothm_0$ in
the finite group
$ 
(\Zee_K/\gothm)^* =
  (\Zee_K/\gothm_0)^* \times \prod_{\gothp|\gothm_\infty} \langle-1\rangle
$
consists of its ordinary residue class modulo $\gothm_0$
and the signs of its images under the real primes $\gothp|\gothm_\infty$.
This group naturally maps onto $P_\gothm/R_\gothm\subset\Cl_\gothm$,
with a kernel
reflecting the fact that generators of principal $\Zee_K$-ideals
are only 
unique up to multiplication by units in $\Zee_K$.

Interpreting both class groups in 2.7 as Galois groups, we see that
all ray class fields contain the Hilbert class field $H=H_1$ from
2.6.2, and that we have an Artin isomorphism
$$
(\Zee_K/\gothm)^*/\im[\Zee_K^*] \isar \Gal(H_\gothm/H)
\leqno{(2.8)}
$$
for their Galois groups over $H$.
By 2.6.1, this generalizes the isomorphism 1.2.

In class field theoretic terms, 
we may specify an abelian extension $K\subset L$ by giving an
admissible modulus $\gothm$ for the extension together with
the corresponding ideal group
$$
A_\gothm = \ker [ I_\gothm \to \Gal(L/K)]
\leqno{(2.9)}
$$
arising as the kernel of the Artin map 2.4.
In this way, we obtain a {\it canonical bijection\/}
between abelian extensions of $K$ inside $\overline K$
and ideal groups $R_\gothm \subset A_\gothm\subset I_\gothm$ of $K$,
provided that one allows for the fact that
the `same' ideal group $A_\gothm$ can be defined modulo different multiples 
$\gothm$
of its {\it conductor\/}, i.e., the conductor of the corresponding extension.
More precisely, we call the ideal groups $A_{\gothm_1}$ and $A_{\gothm_2}$
{\it equivalent\/} if they satisfy
$A_{\gothm_1}\cap I_\gothm=A_{\gothm_2}\cap I_\gothm$ for
some common multiple $\gothm$ of $\gothm_1$ and $\gothm_2$.

Both from a theoretical and an algorithmic point of view, 2.5 provides an
immediate description of the ideal group corresponding to $L$
as the {\it norm group\/} $A_\gothm  = N_{L/K}(I_{\gothm\Zee_L})\cdot R_\gothm$
as soon as we are able to find an admissible modulus $\gothm$ for $L$.
In the reverse direction, finding the {\it class field\/} $L$ corresponding to
an ideal group $A_\gothm$ is much harder.
Exhibiting practical algorithms to do so is the principal task of computational
class field theory, and the topic of this paper.
Already in the case of the Hilbert class field $H$ of $K$ from 2.6.2,
we know no `canonical' generator of $H$, and the problem is non-trivial.

\head
3. Local aspects: ideles
\endhead

\noindent
Over $K=\Que$, all abelian Galois groups
are described as quotients of the groups $(\Zee/m\Zee)^*$
for some modulus $m\in\Zee_{\ge1}$.
One may avoid the ubiquitous choice of moduli that arises 
when dealing with abelian fields
by combining the Artin isomorphisms 1.2 at all `finite levels' $m$
into a single {\it profinite\/} Artin isomorphism
$$
\lim_{\Leftarrow \atop m} (\Zee/m\Zee)^* =
\widehat\Zee^* 
\isar
\Gal(\Que_\ab/\Que)
\leqno{(3.1)}
$$
between the unit group $\widehat\Zee^*$
of the profinite completion $\widehat\Zee$ of $\Zee$
and the absolute abelian Galois group of $\Que$.
The group $\widehat\Zee^*$ splits as a product $\prod_p\Zee_p^*$
by the Chinese remainder theorem, and $\Que_\ab$ is obtained 
correspondingly as a compositum of the fields $\Que(\zeta_{p^\infty})$
generated by the $p$-power roots of unity.
The automorphism corresponding to $u=(u_p)_p\in \widehat\Zee^*$ acts as
$\zeta\mapsto \zeta^{u_p}$ on $p$-power roots of unity.
Note that the component group $\Zee_p^*\subset \widehat\Zee^*$
maps to the inertia group at $p$ in any finite quotient $\Gal(L/\Que)$
of $\Gal(\Que_\ab/\Que)$.

For arbitrary number fields $K$, one can take the projective limit in 2.7 
over all moduli and describe $\Gal(K_\ab/K)$ by an exact sequence
$$
\textstyle
1\tto \Zee_K^*\tto \widehat\Zee_K^*
  \times \prod_{\gothp\text{ real}}\langle-1\rangle
\mapright^{\psi_K} \Gal(K_\ab/K)\tto \Cl_K\tto 1,
\leqno{(3.2)}
$$
which deals with the finite primes occurring in
$\widehat\Zee_K^*=\prod_{\gothp\text{ finite}} U_\gothp$ and the infinite primes
in a somewhat asymmetric way.
Here $\psi_K$ maps the element $-1$ at a real prime $\gothp$ to the 
complex conjugation at the extensions of $\gothp$.
The image of $\psi_K$ is the Galois group $\Gal(K_\ab/H)$ over
the Hilbert class field $H$,
which is of finite index $h_K=\#\Cl_K$ in $\Gal(K_\ab/K)$.
For an abelian extension $L$ of $K$ containing $H$,
the image of the component group $U_\gothp\subset \widehat\Zee_K^*$
in $\Gal(L/H)$ is again the inertia group at $\gothp$ in $\Gal(L/K)$.
As $H$ is totally unramified over $K$, the same is true in case 
$L$ does not contain $H$: the inertia groups for $\gothp$ in
$\Gal(LH/K)$ and $\Gal(L/K)$ are isomorphic under the restriction map.

A more elegant description of $\Gal(K_\ab/K)$ than that provided by
the sequence 3.2 is obtained if one treats
all primes of $K$ in a uniform way, and redefines the Artin map $\psi_K$,
as we will do in 3.7, using the {\it idele group}
$$
\AA_K^*= \textstyle\prod^\prime_\gothp K_\gothp^*
     =\{(x_\gothp)_\gothp : 
        x_\gothp\in U_\gothp \text{ for almost all $\gothp$}\}
$$
of $K$. 
This group [\ST, Section 14]
consists of those elements in the Cartesian product of the multiplicative
groups $K_\gothp^*$ at {\it all\/} completions $K_\gothp^*$ of $K$ that
have their $\gothp$-component in the local unit group $U_\gothp$
for almost all $\gothp$.
Here $U_\gothp$ is, as before, the unit group of the valuation ring
at $\gothp$ in case $\gothp$ is a finite prime of $K$. 
For infinite primes $\gothp$, the choice of $U_\gothp$ is irrelevant
as there are only finitely many such $\gothp$.
We take $U_\gothp=K_\gothp^*$, and write $U_\infty$ to denote
$\prod_{\gothp\text{ infinite}} K_\gothp^*=K\otimes_\Que \RR$.
Note that we have $\prod_{\gothp\text{ finite}} U_\gothp^*=\widehat\Zee_K^*$.

The topology on $\AA_K^*$ is the {\it restricted\/} product topology:
elements are close if they are $\gothp$-adically close at finitely
many $\gothp$
{\it and\/} have a quotient in $U_\gothp$ for all other~$\gothp$.
With this topology,
$K^*$ embeds diagonally into $\AA_K^*$ as a discrete subgroup.
As the notation suggests, $\AA_K^*$ is the unit group of the 
{\it adele ring\/} $\AA_K=\prod^\prime_\gothp K_\gothp$, the subring
of $\prod_\gothp K_\gothp$ consisting of elements having integral components
for almost all $\gothp$.

To any idele $(x_\gothp)_\gothp$, we can associate an ideal
$x\Zee_K=\prod_{\gothp\text{ finite}} \gothp^{\ord_\gothp(x_\gothp)}$,
and this makes the group $I_K$ of fractional $\Zee_K$-ideals into
a quotient of $\AA_K^*$.
For a global element $x\in K^*\subset \AA_K^*$, the ideal $x\Zee_K$ 
is the principal $\Zee_K$-ideal generated by $x$, so we have
an exact sequence
$$
\textstyle
1\tto \Zee_K^*\tto \widehat\Zee_K^* \times U_\infty \tto
\AA_K^*/K^* \tto \Cl_K\tto 1
\leqno{(3.3)}
$$
that describes the {\it idele class group\/} $\AA_K^*/K^*$ of $K$
in a way reminiscent of 3.2.

In order to obtain $\Gal(K_\ab/K)$ as a quotient of $\AA_K^*/K^*$,
we show that the ray class groups
$\Cl_\gothm$ defined in the previous section
are natural quotients of $\AA_K^*/K^*$.
To do so, we associate to a modulus $\gothm=\gothm_0\gothm_\infty$ of $K$
an open subgroup $W_\gothm\subset \AA_K^*$ in the following way.
Write $\gothm=\prod_\gothp \gothp^{n(\gothp)}$ as a formal product,
with $n(\gothp)=\ord_\gothp(\gothm_0)$ for finite $\gothp$, and
$n(\gothp)\in\{0, 1\}$ to indicate the  infinite $\gothp$ in $\gothm_\infty$.
Now put
$$
W_\gothm =\textstyle\prod_\gothp U_\gothp^{(n(\gothp))}
$$
for subgroups $U_\gothp^{(k)}\subset K_\gothp^*$ that are
defined by
$$
U_\gothp^{(k)}=
\cases
U_\gothp&\text{if $k=0$;}\cr
1+\gothp^k&\text{if $\gothp$ is finite and $k>0$;}\cr
U_\gothp^+\subset U_\gothp=\RR^*&\text{if $\gothp$ is real and $k=1$.}\cr
\endcases
$$
Here we write $U_\gothp^+$ for real $\gothp$ to denote the subgroup
of positive elements in $U_\gothp$.
As $\CC^*$ and $\RR_{>0}^*$ have no proper open subgroups,
one sees from the definition of
the restricted product topology on $\AA_K^*$ that a subgroup
$H\subset \AA_K^*$ is open if and only if it contains 
$W_\gothm$ for some modulus $\gothm$.
\proclaim{3.4. Lemma}
For every modulus $\gothm=\prod_\gothp \gothp^{n(\gothp)}$ of $K$, there is an isomorphism
$$\AA_K^*/K^*W_\gothm \isar \Cl_\gothm$$
that maps 
$(x_\gothp)_\gothp$ to the class of 
$\prod_{\gothp\text{ finite}} \gothp^{\ord_\gothp(yx_\gothp)}$.
Here $y\in K^*$ is a global element satisfying 
$yx_\gothp\in U_\gothp^{n(\gothp)}$ for all $\gothp|\gothm$.
\endproclaim\noindent
{\bf Proof.}
Note first that the global element $y$ required in the definition
exists by the approximation theorem. The precise choice of $y$ is irrelevant, as
for any two elements $y, y'$ satisfying the requirement, we have
$y/y'\congr 1\mmod \gothm$.
We obtain a homomorphism $\AA_K^*\to\Cl_\gothm$ that is surjective as it maps
a prime element $\pi_\gothp$ at a finite prime $\gothp\nmid\gothm$ to the
class of $\gothp$.
Its kernel consists of the ideles that can be multiplied into $W_\gothm$ by
a global element $y\in K^*$.\hfill$\square$
\medskip\noindent
If $\gothm$ is an admissible
modulus for the finite abelian extension $K\subset L$, we can
compose the isomorphism in 3.4 with the Artin map 2.4 for $K\subset L$
to obtain an idelic Artin map
$$
\widehat\psi_{L/K}:\qquad 
\AA_K^*/K^* \tto \Gal(L/K)
\leqno{(3.5)}
$$
that no longer refers to the choice of a modulus $\gothm$.
This map, which exists as a corollary of 2.3,
is a continuous surjection that maps the class of a prime element
$\pi_\gothp\in K_\gothp^*\subset \AA_K^*$
to the Frobenius automorphism $\Fr_\gothp\in \Gal(L/K)$ whenever $\gothp$
is finite and unramified in $K\subset L$.

For a finite extension $L$ of $K$, the
adele ring $\AA_L$  is obtained from $\AA_K$ by a base change $K\subset L$, 
so we have a norm map $N_{L/K}:\AA_L\to\AA_K$ that
maps $\AA_L^*$ to $ \AA_K^*$ and restricts to the field norm on
$L^*\subset \AA_L^*$.
As it induces the ideal norm $I_L\to I_K$ on the quotient $I_L$ of $ \AA_K^*$, 
one deduces that the kernel of 3.5
equals $(K^*\cdot N_{L/K} [\AA^*_L]) \mod K^*$, and that we have
isomorphisms
$$
\AA_K^*/K^*N_{L/K} [\AA^*_L] \iso I_\gothm/A_\gothm \isar \Gal(L/K),
\leqno{(3.6)}
$$
with $A_\gothm$ the ideal group modulo $\gothm$ that corresponds to $L$
in the sense of 2.9.
Taking the limit in 3.5 over all finite abelian extensions $K\subset L$
inside $\overline K$, one obtains
the idelic Artin map
$$
\psi_K:\quad \AA_K^*/K^*\tto G_K^\ab=\Gal(K_\ab/K).
\leqno{(3.7)}
$$
This is a continuous surjection 
that is uniquely determined by the property
that the $\psi_K$-image of the class of a prime element
$\pi_\gothp\in K_\gothp^*\subset \AA_K^*$ maps to the Frobenius
automorphism $\Fr_\gothp\in \Gal(L/K)$ for every finite
abelian extension $K\subset L$ in which $\gothp$ is unramified.
It exhibits all abelian Galois groups over $K$
as a quotient of the idele class group $\AA_K^*/K^*$
of $K$.

The kernel of the Artin map 3.7 is the connected component of the unit
element in $\AA_K^*/K^*$.
In the idelic formulation, 
the finite abelian extensions of $K$ inside $\overline K$
correspond bijectively to the open subgroups of $\AA_K^*/K^*$
under the map 
$$
L\longmapsto \psi^{-1}_K[\Gal(K_\ab/L)]
=(K^*\cdot N_{L/K} [\AA^*_L]) \mod K^*.
$$
In this formulation,
computational class field theory amounts to generating, for any given open
subgroup of $\AA_K^*/K^*$, the abelian extension $K\subset L$ 
corresponding to it.
\medskip\noindent
{\bf 3.8. Example.}
Before continuing, let us see what the idelic reformulation of 3.1 
comes down to for $K=\Que$.
Every idele $x=((x_p)_p, x_\infty)\in \AA_\Que^*$ can uniquely be written as
the product of the rational number
$\sign(x_\infty)\cdot \prod_p p^{\ord_p(x_p)}\in\Que^*$ and 
a `unit idele'
$u_x\in\RR_{>0}\times \prod_p\Zee_p^*=\RR_{>0}\times \widehat\Zee^*$.
In this way,
the Artin map 3.7 becomes a continuous surjection
$$
\psi_\Que:\qquad
\AA_\Que^*/\Que^* \iso \widehat\Zee^*\times \RR_{>0}
\tto \Gal(\Que_{\ab}/\Que).
$$
Its kernel is the connected component
$\RR_{>0}\times \{1\}$ of the unit element in $\AA_\Que^*/\Que^*$.
Comparison with 3.1 leads to a commutative diagram of isomorphisms
$$
\matrix
\widehat \Zee^* &\mapright^{-1}&\widehat \Zee^* \cr
\mapdown^{\text{can}}_\wr& &\mapdown^{\text{(3.1)}}_\wr\cr
\AA_\Que^*/(\Que^*\cdot\RR_{>0})&\mapright^{\sim}& \Gal(\Que_\ab/\Que),\cr
\endmatrix
\leqno(3.9)
$$
in which the upper horizontal map is {\it not\/} the identity.
To see this, note that the class of the prime element
$\ell\in\Que_\ell^*\subset \AA_\Que^*$ in 
$\AA_\Que^*/(\Que^*\cdot\RR_{>0})$ is represented by the idele
$x=(x_p)_p\in \widehat\Zee^*$ having components
$x_p=\ell^{-1}$ for $p\ne \ell$ and $x_\ell=1$.
This idele maps to the Frobenius of $\ell$, which
raises roots of unity of order coprime to $\ell$ to their $\ell$-th power.
As $x$ is in all $W_\gothm$ for all conductors $\gothm=\ell^k$, it
fixes $\ell$-power roots of unity.
Thus, the upper isomorphism ``$-1$''
is {\it inversion\/} on $\widehat \Zee^*$.
\endexample
\medskip\noindent
Even though the idelic and the ideal group quotients on the left hand side
of the arrow in 3.6 are the `same' finite group, it is 
the idelic quotient that neatly encodes
information at the {\it ramifying\/} primes $\gothp|\gothm$,
which seem `absent' in the other group.
More precisely, we have for all primes $\gothp$ an injective map
$K_\gothp^*\to \AA_K^*/K^*$ that can be composed with
3.7 to obtain a {\it local\/} Artin map $\psi_{K_\gothp}: K_\gothp^*\to\Gal(L/K)$
at every prime $\gothp$ of $K$.
If $\gothp$ is finite and unramified in $K\subset L$, we have
$U_\gothp\subset \ker \psi_{K_\gothp}$ and an induced isomorphism
of finite cyclic groups
$$
K_\gothp^*/\langle \pi_\gothp^{f_\gothp}\rangle U_\gothp
= K_\gothp^*/ N_{L_\gothq/K_\gothp}[L_\gothq^*]\isar 
\langle \Fr_\gothp\rangle = \Gal(L_\gothq/K_\gothp),
$$
since $\Fr_\gothp$ generates the decomposition group of $\gothp$ in
$\Gal(L/K)$, which
may be identified with the Galois group of the local extension
$K_\gothp\subset L_\gothq$ at a prime $\gothq|\gothp$ in~$L$.
It is a non-trivial fact that 3.5 induces for {\it all\/} primes $\gothp$
of $K$,
also the ramifying and the infinite primes, a {\it local Artin isomorphism\/}
$$
\psi_{L_\gothq/K_\gothp}:
K_\gothp^*/ N_{L_\gothq/K_\gothp}[L_\gothq^*]\isar
\Gal(L_\gothq/K_\gothp).
\leqno(3.10)
$$
In view of our observation after 3.2, it maps
$U_\gothp/N_{L_\gothq/K_\gothp}[U_\gothq]$ for finite $\gothp$
isomorphically onto the inertia group of $\gothp$.

We can use 3.10 to {\it locally\/} compute the exponent $n(\gothp)$ to
which $\gothp$ occurs in the conductor of $K\subset L$:
it is the smallest non-negative integer $k$ for which we have
$U_\gothp^{(k)}\subset N_{L_\gothq/K_\gothp}[L_\gothq^*]$.
For unramified primes $\gothp$ we obtain $n(\gothp)=0$, as the
local norm is surjective on the unit groups. 
For tamely ramified primes we have $n(\gothp)=1$, and
for wildly ramified primes $\gothp$ the exponent $n(\gothp)$
may be found by a local computation.
In many cases it is sufficient to use an upper bound
coming from the fact that every $d$-th power in $K_\gothp^*$
is a norm from $L_\gothq$, with
$d$ the degree of $K_\gothp\subset L_\gothq$ (or even $K\subset L$).
Using Hensel's lemma [\BW], one then finds
$$
n(\gothp)\le e(\gothp/p)\Bigl( {1\over p-1}+ \ord_p (e_\gothp)\Bigr)+1,
\leqno{(3.11)}
$$
where $e(\gothp/p)$ is the absolute ramification index of $\gothp$ over the underlying
rational prime $p$, and $e_\gothp$ the ramification index of $\gothp$ in
$K\subset L$.
Note that $e_\gothp$ is independent of the choice of an extension prime 
as $K\subset L$ is Galois.

\head
4. Computing class fields: preparations
\endhead

\noindent
Our fundamental problem is the computation of the class field $L$ 
that corresponds to a given ideal group $A_\gothm$ of $K$ in the sense of~2.5.
One may `give' $A_\gothm$ by specifying $\gothm$ and a list
of ideals for which the classes in the ray class group $\Cl_\gothm$
generate~$A_\gothm$.
The first step in computing $L$ is the computation of the 
group $I_\gothm/A_\gothm$ that will give us control of the Artin
isomorphism $I_\gothm/A_\gothm\isar\Gal(L/K)$.
As linear algebra over $\Zee$ provides us with good algorithms [\COB, Section 4.1]
to deal with finite or even finitely generated abelian groups,
this essentially reduces to computing the finite group
$\Cl_\gothm$ of which $I_\gothm/A_\gothm$ is quotient.

For the computation of the ray class group $\Cl_\gothm$
modulo $\gothm=\gothm_0\cdot\gothm_\infty$,
one computes, in line with [\SC],
the three other groups in the exact sequence 2.7 in which
it occurs, and the maps between them.
The class group $\Cl_K$ and the unit group $\Zee_K^*$ in 2.7
can be computed using the algorithm described in [\ST, Section 12],
which factors {\it smooth\/} elements of $\Zee_K$ over a {\it factor base\/}.
As this takes exponential time as a function of the base field $K$,
it can only be done for moderately sized $K$.
For the group $(\Zee_K/\gothm_0)^*$, one uses
the Chinese remainder theorem to decompose it into a product of 
local multiplicative groups the form $(\Zee_K/\gothp^k)^*$.
Here we need to assume that we are able
to factor $\gothm_0$, but this is a safe assumption as
we are unlikely to encounter extensions for which we cannot
even factor the conductor.
The group $(\Zee_K/\gothp^k)^*$ is a product of
the cyclic group $k_\gothp^*=(\Zee_K/\gothp)^*$ and the subgroup
$(1+\gothp)/(1+\gothp^k)$, the structure of which can be found inductively
using the standard isomorphisms
$(1+\gothp^a)/(1+\gothp^{a+1})\iso k_\gothp$ and, more efficiently,
$$(1+\gothp^a)/(1+\gothp^{2a})\isar \gothp^a/\gothp^{2a}$$
between multiplicative and additive quotients.
In many cases, the result can be obtained in one stroke using the
$\gothp$-adic logarithm [\COB, Section 4.2.2].
Finding $\Cl_\gothm$ from the other groups in 2.7
is now a standard application of linear algebra over $\Zee$.
The quotient $I_\gothm/A_\gothm$ gives us an explicit description of the
Galois group $\Gal(L/K)$ in terms of Artin symbols of $\Zee_K$-ideals.

For the ideal group $A_\gothm$, we next compute its
conductor $\gothf$, which may be a proper divisor of $\gothm$.
This comes down to checking whether we have 
$A_\gothm\supset I_\gothm\cap R_\gothn$ for some modulus $\gothn|\gothm$.
Even in the case $A_\gothm=R_\gothm$ the conductor can be
smaller than $\gothm$, as the trivial isomorphism
$(\Zee/6\Zee)^*\isar (\Zee/3\Zee)^*$
of ray class groups over $K=\Que$ shows.
The conductor $\gothf$ obtained, which is the same as the conductor
$\gothf_{L/K}$ of the corresponding extension,
is exactly divisible by the primes that ramify in $K\subset L$.
In particular, we know the signature of $L$ from the real primes
dividing $\gothf$.
With some extra effort, one can even compute the discriminant
$\Delta_{L/K}$ using Hasse's
{\it F\"uhrerdiskriminantenproduktformel\/}
$$
\Delta_{L/K}= \prod_{\chi: I_\gothm/A_\gothm\to\CC^*} \gothf(\chi)_0.
\leqno(4.2)
$$
Here $\chi$ ranges over the characters of the finite group
$I_\gothm/A_\gothm\iso\Gal(L/K)$, and $\gothf(\chi)_0$ denotes the 
finite part of the conductor $\gothf(\chi)$ of the ideal group 
$A_\chi$ modulo $\gothm$ satisfying $A_\chi/A_\gothm=\ker\chi$.
All these quantities can be computed by the standard algorithms
for finite abelian groups.
\medskip\noindent
{\bf 4.3. Example.}
If $K\subset L$ is cyclic of prime degree $\ell$, we have
a trivial character of conductor $(1)$ and
$\ell-1$ characters of conductor $\gothf_{L/K}$, so 
4.2 reduces to 
$$\Delta_{L/K}=(\gothf_{L/K})_0^{\ell-1}.$$
In particular, we see that the discriminant of a quadratic extension 
$K\subset L$ is not only for $K=\Que$, but generally
equal to the finite part of the conductor of the extension.
\endexample
\medskip\noindent
Having at our disposal the Galois group $\Gal(L/K)$,
the discriminant $\Delta_{L/K}$ and the Artin isomorphism
$I_\gothm/A_\gothm\isar \Gal(L/K)$ describing
the splitting behavior of the primes in $K\subset L$,
we proceed with the computation of a generator for 
$L$ over $K$, i.e.,
an irreducible polynomial in $K[X]$ with the property that its roots in
$\overline K$ generate $L$.

As the computation of class fields is not an easy computation,
it is often desirable to decompose $\Gal(L/K)$ as a product
$\prod_i \Gal(L_i/K)$ of Galois groups $\Gal(L_i/K)$
and to realize $L$ as a compositum of extensions $L_i$ that
are computed separately.
This way one can work with extensions $L/K$ that are
cyclic of prime power order, or at least of prime power exponent.
The necessary reduction of the global class field theoretic data
for $L/K$ to those for each of the $L_i$ is only a short computation
involving finite abelian groups.

\head
5. Class fields as Kummer extensions
\endhead

\noindent
Let $K$ be {\it any\/}
field containing a primitive $n$-th root of unity $\zeta_n$,
and $K\subset L$ an abelian extension of {\it exponent\/} dividing $n$.
In this situation, Kummer theory [\LA, Chapter VIII, \S6--8] tells us that $L$
can be obtained by adjoining to $K$ the $n$-th roots
of certain elements of $K$.
More precisely, let $W_L=K^*\cap {L^*}^n$ be the subgroup of $K^*$
of elements that have an $n$-th root in~$L$.
Then we have $L=K(\root n\of {W_L})$, and there is the canonical
{\it Kummer pairing\/}
$$
\eqalign{
\Gal(L/K)\times W_L/{K^*}^n &\longrightarrow \langle \zeta_n\rangle\cr
(\sigma, w) \qquad        &\longmapsto\langle\sigma, w\rangle=
 (w^{1/n})^{\sigma-1}={\sigma(\root n\of w)\over \root n\of w}.\cr
}
\leqno{(5.1)}
$$
By {\it canonical\/}, we mean that the natural action of an
automorphism $\tau\in\Aut(\overline K)$ on the pairing for $K\subset L$
yields the Kummer pairing for $\tau K\subset \tau L$:
$$
\langle\tau\sigma\tau^{-1}, \tau w\rangle=
                 \langle\sigma, w\rangle^\tau.
\leqno(5.2)
$$
The Kummer pairing is {\it perfect\/}, i.e., it induces an isomorphism
$$
W_L/{K^*}^n\isar \Hom(\Gal(L/K),\CC^*).
\leqno(5.3)$$ 
In the case where $\Gal(L/K)$ is cyclic of order $n$, this
means that we have $L=K(\root n\of\alpha)$ 
and $W_L=K^*\cap {L^*}^n=\langle\alpha\rangle\cdot {K^*}^n$
for some $\alpha\in K$.
If $\root n\of\beta$ also generates $L$ over $K$, then $\alpha$
and $\beta$ are powers of each other modulo $n$-th powers.

We will apply Kummer theory to generate the class fields of a
number field $K$.
Thus, let $L$ is the class field of $K$ from Section 4
that is to be computed.
Suppose that we have computed a `small' modulus $\gothf$ for $L$ that
is only divisible by the ramifying primes,
e.g. the conductor $\gothf_{L/K}$, and an ideal group
$A_\gothf$ for $L$ by the methods of Section 4.
With this information, we control the Galois group of our extension
via the Artin isomorphism $I_\gothf/A_\gothf\isar \Gal(L/K)$.
Let $n$ be the exponent of $\Gal(L/K)$.
Then we can directly apply Kummer theory if $K$ contains the
required $n$-th roots of unity; if not, we need to pass to a 
cyclotomic extension of $K$ first.
This leads to a natural case distinction.
\medskip\noindent
{\bf Case 1:} $K$ contains a primitive $n$-th root of unity $\zeta_n$.
\smallskip\noindent
Under the restrictive assumption that $K$ contains $\zeta_n$,
the class field $L$ is a Kummer extension of $K$, and
generating $L=K(\root n\of {W_L})$ comes down to finding generators
for the finite group $W_L/{K^*}^n$.
We first compute a {\it finite\/} group containing $W_L/{K^*}^n$.
This reduction is a familiar ingredient from the
{\it proofs\/} of class field theory [\AT, \CF].
\proclaim{5.4. Lemma}
Let $K\subset L$ be finite abelian of exponent $n$, and assume $\zeta_n\in K$.
Suppose $S$ is a finite set of primes of $K$ containing the infinite primes
such that
\itemitem{$1.$}
$K\subset L$ is unramified outside $S$;
\itemitem{$2.$}
$\Cl_K/\Cl_K^n$ is generated by the classes of the finite primes in $S$.
\vskip0cm\noindent
Then the image of the group $U_S$ of $S$-units in $K^*/{K^*}^n$ is
finite of order $n^{\#S}$, and it contains the group $W_L/{K^*}^n$ from
$(5.1)$.
\endproclaim\noindent
The first condition in 5.4 means that $S$ contains all the primes that
divide our small modulus $\gothf$.
The second condition is automatic if the class number of $K$ is prime to $n$,
and it is implied by the first if the classes of the ramifying primes 
generate $\Cl_K/\Cl_K^n$.
Note that any set of elements of $\Cl_K$ generating $\Cl_K/\Cl_K^n$ actually
generates the full `$n$-part' of the class group, i.e., the
product of the $p$-Sylow subgroups of $\Cl_K$ at the primes $p|n$.
In general, there is a lot of freedom in the choice of
primes in~$S$ outside $\gothf$.
One tries to have $S$ `small' in order to minimize the size $n^{\#S}$
of the group $(U_S\cdot {K^*}^n)/{K^*}^n$ containing $W_L/{K^*}^n$.
\medskip\noindent
{\bf Proof of 5.4.}
By the Dirichlet unit theorem [\ST, Theorem 10.9], the group $U_S$ of $S$-units
of $K$ is isomorphic to $\mu_K\times \Zee^{\#S-1}$.
As $\mu_K$ contains $\zeta_n$, the image
$(U_S\cdot {K^*}^n)/{K^*}^n\iso U_S/U_S^n$ of $U_S$ in
$K^*/{K^*}^n$ is finite of order $n^{\#S}$.

In order to show that $(U_S\cdot {K^*}^n)/{K^*}^n$ contains
$W_L/{K^*}^n$, pick any $\alpha\in W_L$.
As $K\subset K(\root n\of\alpha)$ is unramified outside $S$, we have
$(\alpha)=\gotha_S\gothb^n$ for some product $\gotha_S$ of prime ideals
in $S$ and $\gothb$ coprime to all finite primes in $S$.
As the primes in $S$ generate the $n$-part of $\Cl_K$, we can write
$\gothb=\gothb_S\gothc$ with $\gothb_S$ a product of prime ideals
in $S$ and $\gothc$ an ideal of which the class in $\Cl_K$ is of order $u$
coprime to $n$.
Now $\alpha^u$ generates an ideal of the form
$(\alpha^u)=\gotha_S'(\gamma^n)$ with $\gotha_S'$
a product of prime ideals in $S$ and $\gamma\in K^*$.
It follows that $\alpha^u\gamma^{-n}\in K^*$ is an $S$-unit, so
$\alpha^u$ and therefore $\alpha$ is contained in $U_S\cdot {K^*}^n$.
\hfill$\square$
\medskip\noindent
In the situation of 5.4, we see that $K\subset L$ is a subextension of
the Kummer extension $K\subset N=K(\root n\of {U_S})$ of degree $n^{\#S}$.
We have to find the subgroup of $U_S/U_S^n$ corresponding to $L$.
This amounts to a computation in linear algebra using the Artin map
and the Kummer pairing.
For ease of exposition, we assume that the set $S$ 
we choose to satisfy 5.4 contains all primes dividing $n$.
This implies that $N$ is the maximal abelian extension of exponent $n$
of $K$ that is unramified outside~$S$.

As we compute $L$ as a subfield of the abelian extension $K\subset N$, we replace
the modulus $\gothf$ of $K\subset L$ by some multiple $\gothm$ that is
an admissible modulus for $K\subset N$.
Clearly $\gothm$ only needs to be divisible by the ramified primes
in $K\subset N$, which are all in $S$.
Wild ramification only occurs at primes $\gothp$ dividing $n$, and
for these primes we can take $\ord_\gothp(\gothm)$ equal
to the bound given by 3.11.
The ideal group modulo $\gothm$ corresponding to $N$ is 
$I_{\gothm}^n\cdot P_{\gothm}$ as $N$ is the maximal exponent $n$ extension
of $K$ of conductor $\gothm$,
and the Artin map for $K\subset N$ is
$$
I_\gothm\tto
I_{\gothm}/(I_{\gothm}^n\cdot P_{\gothm})=\Cl_{\gothm}/\Cl_{\gothm}^n
\isar \Gal(N/K).
\leqno{(5.5)}
$$
The induced map $I_\gothm\to \Gal(L/K)$ is the Artin map
for $K\subset L$, which has the ideal group $A_\gothm$ corresponding
to $L$ as its kernel.
Let $\Sigma_L\subset I_{\gothm}$ be a finite set of ideals of which
the classes generate the $\Zee/n\Zee$-module
$A_\gothm/(I_{\gothm}^n\cdot P_{\gothm})\iso \Gal(N/L)$.
We then have to determine the subgroup $V_L\subset U_S$ consisting of those
$S$-units $v\in U_S$ that have the property that
$\root n\of v$ is left invariant by the Artin symbols of all ideals in
$\Sigma_L$, since the class field we are after is $L=K(\root n\of {V_L})$.

We are here in a situation to apply linear algebra over $\Zee/n\Zee$,
as the Kummer pairing 5.1 tells us that the action of the Artin symbols 
$\Art_{N/K}(\gotha)$ of the ideals $\gotha\in I_\gothm$ on the $n$-th
roots of the $S$-units is described by the pairing of $\Zee/n\Zee$-modules
$$
\eqalign{
I_{\gothm}/I_{\gothm}^n\times U_S/U_S^n &
        \longrightarrow \langle\zeta_n\rangle\cr
(\gotha, u)&\longmapsto \langle\Art_{N/K}(\gotha),u\rangle=
                       (u^{1/n})^{\Art_{N/K}(\gotha)-1}.
}\leqno{(5.6)}
$$
Making this computationally explicit amounts to computing the pairing for
some choice of basis elements of the three modules involved.

For $\langle\zeta_n\rangle$ we have the obvious $\Zee/n\Zee$-generator
$\zeta_n$,
and $I_{\gothm}/I_{\gothm}^n$ is a free $\Zee/n\Zee$-module generated 
by the primes $\gothp\notin S$.
If $K$ is of moderate degree, the general algorithm
[\ST, Section 12] for computing units and class groups can be used to
compute generators for $U_S$, which then form a $\Zee/n\Zee$-basis
for $U_S/U_S^n$.
In fact, finding $s-1=\#S-1$ independent units in $U_S$ that generate
a subgroup of index coprime to $n$ is enough: together with a 
root of unity generating $\mu_K$, these will generate $U_S/U_S^n$.
This is somewhat easier than finding actual generators for $U_S$,
as maximality modulo $n$-th powers is not difficult to establish
for a subgroup $U\subset U_S$ having the right rank $r=\text{rank}_{\Zee/n\Zee}(U)$. 
Indeed, each reduction modulo a small prime
$\gothp\notin S$ provides a character $U\subset U_S\to k_\gothp^*/(k_\gothp^*)^n\iso
\langle\zeta_n\rangle$, the $n$-th power residue symbol at $\gothp$.
By finding $r$ independent characters, one shows that
the intersection of their kernels equals $U^n=U\cap U_S^n$.

For a prime $\gothp\notin S$ and $u\in U_S$, the definitions of
the Kummer pairing and the Frobenius automorphism yield
$$
\langle \Fr_\gothp, u\rangle=
(u^{1/n})^{\Fr_\gothp-1}\congr u^{(N\gothp-1)/n} \in k_\gothp^*,
$$
with $N\gothp=\#k_\gothp$ the absolute norm of $\gothp$.
Thus $\langle \Fr_\gothp, u\rangle$ is simply the power of $\zeta_n$
that is congruent to $u^{(N\gothp-1)/ n}\in k_\gothp^*$.
Even when $\gothp$ is large, this
is not an expensive discrete logarithm problem in $k_\gothp^*$,
since in practice the exponent $n\le [L:K]$ is small:
one can simply check all powers of $\overline \zeta_n\in k_\gothp^*$.
As $\langle\zeta_n\rangle$ reduces injectively modulo 
primes $\gothp\nmid n$, the $n$-root of unity
$\langle \Fr_\gothp, u\rangle$ can be recovered
from its value in $k_\gothp^*$.

{}From the values $\langle \Fr_\gothp, u\rangle$, we compute all symbols
$\langle\Art_{N/K}(\gotha),u\rangle$ by linearity.
It is now a standard computation in linear algebra to find 
generators for the subgroup $V_L/U_S^n\subset U_S/U_S^n$ that is
annihilated by the ideals $\gotha\in\Sigma_L$ under the pairing 5.6.
This yields explicit generators for the Kummer extension
$L=K(\root n\of {V_L})$, and concludes the computation of $L$ in the
case where $K$ contains $\zeta_n$, with $n$ the exponent of $\Gal(L/K)$.
\medskip\noindent
{\bf Case 2:} $K$ does not contain $\zeta_n$.
\smallskip\noindent
In this case $L$ is not a Kummer extension of $K$,
but $L'=L(\zeta_n)$ is a Kummer extension of $K'=K(\zeta_n)$.
$$
\def\cross#1#2{\setbox0\hbox{$#1$}%
  \hbox to\wd0{\hss\hbox{$#2$}\hss}\llap{\unhbox0}}
\varrowlength=15pt
\sarrowlength=16pt
\gridcommdiag{
&&L'\rlap{$=L(\zeta_n)$}&&&\cr
&\sline(1,1)&&{\sline(1,-1)}&&\cr
L&&&&\sline(1,-1)&\cr
&\sline(1,-1)&&&&K'\rlap{$=K(\zeta_n)$}\cr
&&\sline(1,-1)&&\sline(1,1)&\cr
&&&L\cap K'&&\cr
&&&\vline&&\cr
&&&K&&\cr
}
\leqno{(5.7)}
$$
In order to find generators of $L'$ over $K'$ by the method of
Case 1, we need to `lift'
the class field theoretic data from $K$ to $K'$ to
describe $L'$ as a class field of $K'$.
Lifting the modulus $\gothf=\gothf_0\gothf_\infty$ for $K\subset L$
is easy: as $K'$ is totally complex,
$\gothf'=\gothf_0\Zee_{K'}$ is admissible for $K'\subset L'$.
{}From the definition of the Frobenius automorphism,
it is immediate that we have a commutative diagram
$$
\matrix
I_{K',\gothf'} &\mapright^{\text{Artin}}& \Gal(L'/K')\cr
\mapdown_{N_{K'/K}}& &\mapdown_{\text{res}}\cr
I_{K,\gothf}&\mapright^{\text{Artin}}& \Gal(L/K).\cr
\endmatrix
$$
As the restriction map on the Galois groups is injective, we see that
the inverse norm image $N_{K'/K}^{-1}A_\gothf\subset I_{K', \gothf'}$
is the ideal group of $K'$ corresponding to the extension $K'\subset L'$.
As $N_{K'/K}^{-1}A_\gothf$ contains $P_{K', \gothf'}$,
computing this inverse image takes place inside the finite group
$\Cl_{K', \gothf'}$, a ray class group for $K'$.

We perform the algorithm from Case 1 for the extension
$K'\subset L'$ in order to find generators of $L'$ over $K'$.
We are then working with (ray) class groups and $S$-units 
in $K'$ rather than in $K$, and $S$ has to satisfy 5.4.2 for
$\Cl_{K'}$.
All this is only feasible if $K'$ is of moderate degree, and this
seriously restricts the values of $n$ one can handle in practice.
Our earlier observation that we may
decompose $I_\gothf/A_\gothf\iso \Gal(L/K)$
into a product of cyclic groups of prime power order and generate $L$
accordingly as a compositum of cyclic extensions of $K$
is particularly relevant in this context, as it reduces
our problem to a number of instances where $K\subset L$
is cyclic of prime power degree.
Current implementations [\FI] deal with prime power values up to $20$.

We further assume for simplicity that we are indeed in the case where
$K\subset L$ is cyclic of prime power degree~$n$, with $K'=K(\zeta_n)\ne K$.
Suppose that, using the algorithm from Case 1,
we have computed a Kummer generator
$\theta\in L'$ for which we have $L'=K'(\theta)=K(\zeta_n, \theta)$ and
$\theta^n=\alpha\in K'$.
We then need to `descend' $\theta$ efficiently
to a generator $\eta$ of $L$ over $K$.
If $n$ is prime, one has $L=K(\eta)$ for the trace
$$
\eta=\Tr_{L'/L}(\theta).
\leqno{(5.8)}
$$
For prime powers this does not work in all cases.
One can however replace $\theta$ by $\theta+k\zeta_n$
for some small integer $k\in\Zee$ to ensure that
$\theta$ generates $L'$ over $K$, and then general field theory tells us that
the coefficients of the irreducible polynomial
$$f_L^\theta=\prod_{\tau\in \Gal(L'/L)} (X-\tau(\theta))\in L[X]
\leqno{(5.9)}
$$
of $\theta$ over $L$ generate $L$ over $K$.
As we took $K\subset L$ to be cyclic of prime power degree,
one of the coefficients
is actually a generator, and in practice the trace works.
In all cases, one needs an explicit description of the action
of the Galois group $\Gal(L'/L)$ on $\theta$ and $\zeta_n$ in order to compute
the trace 5.8, and possibly other coefficients of $f_L^\theta$ in 5.9.
Finally, if we have $L=K(\eta)$, we need the action of $\Gal(L/K)$ on
$\eta$ in order to write down the generating polynomial 
$$f^\eta_K=\prod_{\sigma\in \Gal(L/K)} (X-\sigma(\eta))\in K[X]$$
for $K\subset L$ that we are after.

As before, the Artin map gives us complete control over the action
of the abelian Galois group $\Gal(L'/K)$ on $L'=K(\zeta_n, \theta)$,
provided that we describe the elements of $\Gal(L'/K)$ as Artin symbols.
We let $\gothm$ be an admissible modulus for $K\subset L'$;
the least common multiple of $\gothf_{L/K}$ and
$n\cdot\prod_{\gothp\text{ real}} \gothp$
is an obvious choice for $\gothm$.
All we need to know is the explicit action of the Frobenius automorphism
$\Fr_\gothp\in\Gal(L'/K)$ of a prime $\gothp\nmid \gothm$ of $K$ on
the generators $\zeta_n$ and $\theta$ of $L'$ over $K$.
Note that $\gothp$ does not divide $n$, and that we may assume that 
$\alpha=\theta^n$ is a unit at $\gothp$.

The cyclotomic action of $\Fr_\gothp\in\Gal(L'/K)$ is given by
$\Fr_\gothp(\zeta_n)=\zeta_n^{N\gothp}$, with $N\gothp=\#k_\gothp$
the absolute norm of $\gothp$.
This provides us with the Galois action on $K'$, and yields
canonical isomorphisms
$$\eqalignno{
\Gal(K'/K)&\iso \im [N_{K/\Que}: I_\gothm\tto (\Zee/n\Zee)^*],\cr
\Gal(K'/(L\cap K'))&\iso \im [N_{K/\Que}: A_\gothm\tto (\Zee/n\Zee)^*].\cr
}$$
%
In order to understand the action of $\Fr_\gothp$ on $\theta=\root n\of\alpha$, we
first observe that $K\subset L'=K'(\theta)$ can only be abelian if $\alpha$ is
in the {\it cyclotomic\/} eigenspace of $(K')^*$ modulo $n$-th powers 
under the action of $\Gal(K'/K)$.
More precisely, applying 5.2 for $K'\subset L'$ with
$\tau= \Fr_\gothp$, we have
$\Fr_\gothp\cdot\sigma\cdot \Fr_\gothp^{-1}= \sigma$ as $\Gal(L'/K)$ is abelian,
and therefore 
$$
\langle\sigma, \Fr_\gothp(\alpha)\rangle=
\langle\sigma, \alpha\rangle^{\Fr_\gothp}=
\langle\sigma, \alpha\rangle^{N\gothp}=\langle\sigma, \alpha^{N\gothp}\rangle
$$
for all $\sigma\in \Gal(L'/K')$.
By 5.3, we conclude that we have 
$$
\Fr_\gothp(\alpha)= \alpha^{N\gothp} \cdot \gamma_\gothp^n
\leqno{(5.10)}
$$
for some element $\gamma_\gothp\in K'$.
As we know how $\Fr_\gothp$ acts on $\alpha\in K'=K(\zeta_n)$, we can
compute $\gamma_\gothp$ by extracting
some $n$-th root of $\Fr_\gothp(\alpha)\alpha^{-N\gothp}$ in $K'$.
The element $\gamma_\gothp$ is only determined up to multiplication
by $n$-th roots of unity by 5.10. 
As we took $\alpha$ to be a unit at $\gothp$,
we have $\gamma_\gothp^n\congr1\mod\gothp$
by definition of the Frobenius automorphism, so 
there is a unique element $\gamma_\gothp\congr1\mod\gothp$ satisfying 5.10.
With this choice of $\gamma_\gothp$, we have
$$
\Fr_\gothp(\theta)= \theta^{N\gothp} \cdot \gamma_\gothp
$$
as the $n$-th powers of both quantities are the same by 5.10, and they are
congruent modulo $\gothp$.
This provides us with the explicit Galois action of $\Fr_\gothp$ on $\theta$
for unramified primes $\gothp$.

The description of the Galois action on $\theta$ and $\zeta_n$ in terms of
Frobenius symbols is all we need.
The Galois group $\Gal(L'/L)\iso \Gal(K'/(L\cap K'))$, which we may identify
with the subgroup $N_{K/\Que}(A_\gothm)$ of $(\Zee/n\Zee)^*$, is either cyclic
or, in case $n$ is a power of 2, generated by 2 elements.
Picking one or two primes $\gothp$ in $A_\gothm$ with norms in suitable 
residue classes
modulo $n$ is all it takes to generate $\Gal(L'/L)$ by Frobenius automorphisms,
and we can use these elements to descend $\theta$ to a generator $\eta$ for 
$L$ over $K$.
We also control the Galois action of $\Gal(L/K)=I_\gothm/A_\gothm$ on 
$\eta$, and this
makes it possible to compute the irreducible polynomial $f^\eta_K$ for the
generator $\eta$ of $L$ over $K$.

\head
6. Class fields arising from complex multiplication
\endhead

\noindent
As we observed in Example 2.6.1, the ray class fields over
the rational number field~$\Que$ are the cyclotomic fields.
For these fields, we have explicit generators over $\Que$ that
arise `naturally' as the values of the analytic function
$q: x\mapsto e^{2\pi i x}$
on the unique archimedean completion $\RR$ of $\Que$.
The function $q$ is periodic modulo the ring of integers
$\Zee\subset \RR$ of~$\Que$, and it induces an isomorphism
$$
\eqalign{
\RR/\Zee &\  \isar T=\{z\in\CC: z\overline z=1\}\subset\CC\cr
     x&\longmapsto q(x)=e^{2\pi i x}\cr
}
\leqno(6.1)
$$
between the quotient group $\RR/\Zee$ and the `circle group' $T$
of complex numbers of absolute value 1.
The Kronecker-Weber theorem 1.1 states that the values of the analytic
function $q$ at the points of the {\it torsion subgroup\/}
$\Que/\Zee\subset\RR/\Zee$ generate the maximal abelian extension
$\Que_\ab$ of $\Que$.
More precisely,
the $q$-values at the $m$-torsion subgroup ${1\over m}\Zee/\Zee$ of
$\RR/\Zee$ generate the $m$-th cyclotomic field $\Que(\zeta_m)$.
Under this parametrization of roots of unity by $\Que/\Zee$,
the Galois action on the $m$-torsion values comes from 
multiplications on ${1\over m}\Zee/\Zee$ by integers $a\in\Zee$ coprime to $m$,
giving rise to the Galois group 
$$\textstyle
\Gal(\Que(\zeta_m)/\Que)=\Aut({1\over m}\Zee/\Zee)=(\Zee/m\Zee)^*
\leqno{(6.2)}
$$
from 1.2.
Taking the projective limit over all $m$, one obtains the identification
of $\Gal(\Que_\ab/\Que)$ with $\Aut(\Que/\Zee)=\widehat\Zee^*$
from 3.1, and we saw in 3.8 that the relation with 
the Artin isomorphism is given by the commutative diagram 3.9.
To stress the analogy with the complex multiplication case, we rewrite 3.9 as
$$
\matrix
\widehat \Zee^* &\mapright^{-1}&
                \widehat \Zee^*\rlap{$=\AA_\Que^*/(\Que^*\cdot\RR_{>0})$}\cr
\mapdown^{\text{can}}_\wr& &\mapdown^{\text{Artin}}_\wr\cr
\Aut(\Que/\Zee)&\mapright^{\sim}& \Gal(\Que_\ab/\Que),\cr
\endmatrix
\leqno(6.3)
$$
where ``$-1$" denotes inversion on $\widehat \Zee^*$.

From now on, we take $K$ to be an imaginary
quadratic field.
Then $K$ has a single archimedean completion $K\to\CC$,
and much of what we said for the analytic function
$q$ on $\RR/\Zee$ has an analogue for the quotient group $\CC/\Zee_K$.
In complete analogy, we will define an analytic function
$f_K: \CC/\Zee_K\to \PP^1(\CC)$ in 6.14 with the property that its 
finite values at the $m$-torsion subgroup ${1\over m}\Zee_K/\Zee_K$ of
$\CC/\Zee_K$ generate the ray class field $H_m$ of $K$
of conductor $m\Zee_K$.
However, to define this {\it elliptic function\/} $f_K$ on 
the complex {\it elliptic curve\/} $\CC/\Zee_K$, we need an 
algebraic description of $\CC/\Zee_K$, which exists over an
{\it extension\/} of $K$ that is usually larger than $K$ itself.
The Hilbert class field $H=H_1$ of~$K$ from 2.6.2
is the smallest extension of $K$ that one can use,
and the torsion values of $f_K$ generate class fields over $H$.
This makes the construction of $H$ itself into an important
preliminary step that does not occur over $\Que$, as $\Que$
is its own Hilbert class field.

In this section, we give the classical algorithms for
constructing the extensions $K\subset H$ and $H\subset H_m$.
The next section provides some theoretical background and 
different views on complex multiplication.
Our final section 8 shows how such views lead to algorithmic
improvements.
\medskip\noindent
Complex multiplication starts with the fundamental observation 
[\SI, Chapter VI] that for every lattice $\Lambda\subset\CC$,
the complex torus $\CC/\Lambda$ admits a meromorphic function, the
Weierstrass $\wp$-function
$$
\textstyle
\wp_\Lambda: z\mapsto z^{-2}+\sum_{\omega\in\Lambda\setminus\{0\}}
[(z-\omega)^{-2}-\omega^{-2}],
$$
that has period lattice $\Lambda$ and is holomorphic except for
double poles at the points of~$\Lambda$.
The corresponding Weierstrass map
$$
\eqalign{
W: \CC/\Lambda&\tto E_\Lambda\subset \PP^2(\CC)\cr
             z&\longmapsto [\wp_\Lambda(z): \wp'_\Lambda(z):1]\cr}
$$
is a complex analytic isomorphism
between the torus $\CC/\Lambda$ and the complex elliptic curve
$E_\Lambda\subset \PP^2(\CC)$ defined by the affine Weierstrass equation
$$
y^2=4x^3-g_2(\Lambda)x-g_3(\Lambda).
$$
The Weierstrass coefficients
$$
\textstyle
g_2(\Lambda)=60\sum\limits_{\omega\in\Lambda\setminus\{0\}} \omega^{-4}
\quad\text{and}\quad
g_3(\Lambda)=140\sum\limits_{\omega\in\Lambda\setminus\{0\}} \omega^{-6}
\leqno{(6.4)}
$$
of $E_\Lambda$ are the {\it Eisenstein series\/} of weight 4 and 6
for the lattice $\Lambda$.
The natural addition on $\CC/\Lambda$ translates into an algebraic
group structure on $E_\Lambda(\CC)$ sometimes referred to as
`chord and tangent addition'.
On the Weierstrass model $E_\Lambda$, 
the point $O=[0:1:0]=W(0\mod\Lambda)$ at infinity is the zero point,
and any line in $\PP^2(\CC)$ intersects the curve $E_\Lambda$
in 3 points, counting multiplicities, that have sum $O$.

All complex analytic maps $\CC/\Lambda_1\to\CC/\Lambda_2$
fixing the zero point are 
multiplications $z\mapsto\lambda z$ with $\lambda\in\CC$
satisfying $\lambda\Lambda_1\subset\Lambda_2$.
These are clearly group homomorphisms, and in the commutative diagram
$$
\matrix
\CC/\Lambda_1&\mapright^{\lambda}& \CC/\Lambda_2\cr
\mapdown^{W_1}_\wr& &\mapdown^{W_2}_\wr\cr
E_{\Lambda_1}&\mapright^{\phi_\lambda}& E_{\Lambda_2}\cr
\endmatrix
\leqno(6.5)
$$
the corresponding maps
$\phi_\lambda: E_{\Lambda_1}\tto E_{\Lambda_2}$
between algebraic curves are known as {\it isogenies\/}.
For $\lambda\ne0$, the isogeny $\phi_\lambda$
is a finite algebraic map of degree
$[\Lambda_2:\lambda\Lambda_1]$,
and $E_{\Lambda_1}$ and $E_{\Lambda_2}$ are isomorphic
as complex algebraic curves
if and only if we have $\lambda\Lambda_1=\Lambda_2$ for some
$\lambda\in\CC$.
The isogenies $E_\Lambda\to E_\Lambda$ form the
{\it endomorphism ring\/}
$$
\End(E_\Lambda)
\iso
\{\lambda\in\CC: \lambda\Lambda\subset\Lambda\}
\leqno(6.6)
$$
of the curve $E_\Lambda$, which we can view as
a discrete subring of $\CC$.
The $\lambda$-value of the analytically defined endomorphism
`multiplication by $\lambda\in\CC$' is reflected
algebraically as a true multiplication by $\lambda$
of the {\it invariant differential\/}
$dx/y$ on $E_\Lambda$ coming from $dz=d(\wp_\Lambda)/\wp_\Lambda'$.
If $\End(E_\Lambda)$ is strictly larger than $\Zee$, it is
a complex quadratic order $\Cal O$ and $E_\Lambda$ is said to have
{\it complex multiplication\/} (CM) by $\Cal O$.

In order to generate the class fields of our imaginary 
quadratic field~$K$,
we employ an elliptic curve $E_\Lambda$ having CM by $\Zee_K$.
Such a curve can be obtained by taking $\Lambda$ equal to $\Zee_K$, or to
a fractional $\Zee_K$-ideal~$\gotha$, but the
Weierstrass coefficients 6.4 for $\Lambda=\gotha$ will not
in general be algebraic.

In order to find an {\it algebraic\/} model for the complex curve
$E_\Lambda$, we scale $\Lambda$ to a {\it homothetic\/} lattice
$\lambda\Lambda$ to obtain a $\CC$-isomorphic model
$$
E_{\lambda\Lambda}:
y^2=4x^3-\lambda^{-4}g_2(\Lambda)x-\lambda^{-6}g_3(\Lambda)
$$
under the Weierstrass map.
The discriminant $\Delta=g_2^3-27g_3^2$
of the Weierstrass polynomial $4x^3-g_2x-g_3$ does not vanish, and
the lattice function
$$\Delta(\Lambda)=g_2(\Lambda)^3-27g_3(\Lambda)^2
$$
is of weight 12:
it satisfies $\Delta(\lambda\Lambda)=\lambda^{-12} \Delta(\Lambda)$.
Thus, the $j$-invariant
$$
j(\Lambda)=1728\, {g_2(\Lambda)^3\over \Delta(\Lambda)}
          =1728\, {g_2(\Lambda)^3\over g_2(\Lambda)^3-27g_3(\Lambda)^2}
\leqno(6.7)
$$
is of weight zero, and an invariant of the homothety
class of $\Lambda$ or, equivalently, the isomorphism class
of the complex elliptic curve $\CC/\Lambda$.
It generates the minimal field of definition over which $\CC/\Lambda$
admits a Weierstrass model.

If $E_\Lambda$ has CM by $\Zee_K$, then $\Lambda$ is homothetic
to some $\Zee_K$-ideal~$\gotha$.
It follows that, up to isomorphism, there are only finitely many 
complex elliptic curves $E_\Lambda$ having CM by $\Zee_K$, one for
each ideal class in $\Cl_K$.
As any automorphism of $\CC$ maps the algebraic curve $E_\Lambda$ to
an elliptic curve with the same endomorphism ring, we find that the
$j$-invariants of the ideal classes of $\Zee_K$ form a set of
$h_K=\#\Cl_K$ distinct {\it algebraic\/} numbers permuted
by the abolute Galois group $G_\Que$ of $\Que$.
This allows us to define the {\it Hilbert class polynomial\/} of $K$
as
$$
\Hil_K(X)=
 \prod_{[\gotha] \in \Cl_K}(X - j(\gotha)) \in \Que[X].
\leqno{(6.8)}
$$
Its importance stems from the following theorem, traditionally
referred to as the {\it first main
theorem of complex multiplication\/}.
\proclaim
{6.9. Theorem}
The Hilbert class field $H$ of $K$ is the splitting field
of the polynomial $\Hil_K(X)$ over $K$. 
This polynomial is irreducible in $K[X]$, and the Galois action of
the Artin symbol $\sigma_\gothc=\Art_{H/K}(\gothc)$
of the ideal class $[\gothc]\in\Cl_K\iso \Gal(H/K)$
on the roots $j(\gotha)$ of $\Hil_K(X)$  is given by
$$
j(\gotha)^{\sigma_\gothc}=j(\gotha\gothc^{-1}).
$$
\endproclaim
\noindent
To compute $\Hil_K(X)$ from its definition 6.8, one compiles a
list of $\Zee_K$-ideal classes in
the style of Gauss, who did this in terms of binary quadratic forms.
Every $\Zee_K$-ideal class $[\gotha]$ has a representative of the form
$\Zee\tau+\Zee$, with $\tau\in K$ a root of some
irreducible polynomial $aX^2+bX+c\in\Zee[X]$  of discriminant
$b^2-4ac=\Delta_{K/\Que}$.
If we take for $\tau$ the root in the complex
upper half plane $\HH$, the {\it orbit\/} of $\tau$ under the
natural action 
$$
\left(
\matrix \alpha&\beta\cr \gamma&\delta\cr \endmatrix
\right)
(z)=
{\alpha z+\beta\over \gamma z+\delta}
$$
of the modular group $\Sl_2(\Zee)$ on $\HH$
is uniquely determined by $[\gotha]\in\Cl_K$.
In this orbit, there is a unique element
$$
\tau_\gotha={-b+\sqrt {b^2-4ac}\over 2a} \in\HH
$$
that lies in the standard fundamental domain 
for the action of $\Sl_2(\Zee)$ on $\HH$ consisting of those
$z\in\HH$ that satisfy the two inequalities $|\text{Re}(z)|\le{1\over 2}$
and $z\overline z\ge1$ and,
in case we have equality in either of them, also $\text{Re}(z)\le0$.
This yields a description
$$
[\gotha]=
\textstyle
[\Zee\cdot {-b+\sqrt {b^2-4ac}\over 2a} +\Zee]
\qquad\longleftrightarrow\qquad
(a,b,c)
$$
of the elements of $\Cl_K$ as {\it reduced\/} integer triples
$(a,b,c)$ of discriminant $b^2-4ac=\Delta_{K/\Que}$.
As we have Re$(\tau_\gotha)=-{b\over 2a}$ and
$\tau_\gotha\overline\tau_\gotha= {c\over a}$,
the reduced integer triples $(a,b,c)$ corresponding to
$\tau_\gotha$ in the fundamental domain for $\Sl_2(\Zee)$
are those satisfying
$$
|b|\le a\le c\qquad\text{and}\qquad b^2-4ac=\Delta_{K/\Que},
$$
where $b$ is non-negative in case we have $|b|=a$ or $a=c$.
For reduced forms, one sees from the inequality
$
\Delta_{K/\Que}=b^2-4ac\le a^2-4a^2=-3a^2
$
that we have $|b|\le a\le \sqrt{|\Delta_{K/\Que}|/3}$,
so the list is indeed finite, and can easily be generated
[\COA, Algorithm 5.3.5].
See [\CX] for the classical interpretation of the triples
$(a,b,c)$ as positive definite integral binary quadratic forms
$aX^2+bXY+cY^2$ of discriminant $b^2-4ac=\Delta_{K/\Que}$.

If we put $j(\tau)=j(\Zee\tau+\Zee)$, 
the $j$-function 6.7 becomes a holomorphic function $j:\HH\to\CC$
invariant under the action of $\Sl_2(\Zee)$.
As it is in particular invariant under $\tau\mapsto \tau+1$,
it can be expressed in various ways in terms of the variable
$q=e^{2\pi i \tau}$ from 6.1.
Among them is the well-known {\it integral\/} Fourier expansion
$$
j(\tau)=j(q)=q^{-1}+744+196884 q+ \ldots \in q^{-1}+\Zee[[q]]
\leqno(6.10)
$$
that explains the normalizing factor 1728 in the definition 6.7 of $j$.
It implies [\LAEF, Chapter~5, \S2] that the roots of
$\Hil_K(X)$ in 6.8 are algebraic
integers, so $\Hil_K(X)$ is a polynomial in $\Zee[X]$ that can
be computed {\it exactly\/} from complex approximations of its roots
that are sufficiently accurate to yield the right hand side of 6.8
in $\CC[X]$ to `one-digit precision'.
For numerical computations of $j(\tau)$, one uses
approximate values of the Dedekind $\eta$-function
$$
\eta(\tau)=q^{1/24}\prod_{n\ge1}(1-q^n)
=q^{1/24}\sum_{n\in\Zee} (-1)^n q^{n(3n-1)/2},
\leqno(6.11)
$$
which has a lacunary Fourier expansion that is better suited for numerical
purposes than 6.10.
From $\eta$-values one computes $\gothf_2(\tau)=\sqrt 2 \eta(2\tau)/\eta(\tau)$,
and finally $j(\tau)$ as
$$
j(\tau)={(\gothf_2^{24}(\tau)+16)^3\over \gothf_2^{24}(\tau)}.
\leqno(6.12)
$$
This finishes the description of the classical algorithm
to compute the Hilbert class field $H$ of $K$.
\medskip\noindent
Having computed the irreducible polynomial $\Hil_K(X)$ of
$j_K=j(\Zee_K)$,
we can write down a Weierstrass model $E_K$ for $\CC/\Zee_K$
over $H=K(j_K)$ (or even over $\Que(j_K)$) and use it to generate the ray 
class field extensions $H\subset H_m$.
Choosing $E_K$ is easy in the special cases $K=\Que(\zeta_3), \Que(i)$, when
one of $g_2=g_2(\Zee_K)$ and $g_3=g_3(\Zee_K)$ vanishes
and the other can be scaled to have any non-zero rational value.
For $K\ne\Que(\zeta_3), \Que(i)$, the number
$\lambda=\sqrt{g_3/g_2}\in\CC^*$ is determined up to sign, and since
we have
$$
c_K=\lambda^{-4}g_2=\lambda^{-6}g_3=g_2^3g_3^{-2}=27 {j_K\over j_K-1728},
$$
the Weierstrass model $y^2=4x^3-c_Kx-c_K$ for $\CC/\lambda\Zee_K$
is defined over $\Que(j_K)\subset H$.
A more classical choice is $\lambda^2=\Delta/(g_2g_3)$,
with $g_2$, $g_3$ and $\Delta$ associated to $\Zee_K$,
giving rise to the model
$$
E_K: y^2= w_K(x), \quad \text{where}\quad
w_K(x)=4x^3-{c_K\over(c_K-27)^2} x- {c_K\over(c_K-27)^3}.
\leqno(6.13)
$$
Any scaled Weierstrass parametrization $W_K:\CC/\Zee_K\isar E_K$ with
$E_K$ defined over $H$ can serve as the imaginary quadratic analogue
of the isomorphism $q: \RR/\Zee\isar T$ in~6.1.
For the model $E_K$ in 6.13,
the $x$-coordinate $\wp_{\lambda\Zee_K}(\lambda z)=
\lambda^{-2}\wp_{\Zee_K}(z)$ of $W_K(z)$ is given by the {\it Weber function}
$$
f_K(z)= {g_2(\Zee_K)g_3(\Zee_K)\over \Delta(\Zee_K)} \wp_{\Zee_K}(z).
\leqno(6.14)
$$
It has `weight 0' in the sense that the right hand side is invariant
under simultaneous scaling $(\Zee_K, z)\to(\lambda \Zee_K, \lambda z)$
of the lattice $\Zee_K$ and the argument $z$ by $\lambda\in\CC^*$.

In the special cases $K=\Que(i)$ and $\Que(\zeta_3)$ that
have $\Zee_K^*$ of order 4 and 6, there are slightly different 
Weber functions $f_K$ that are not the $x$-coordinates on a Weierstrass
model for $\CC/\Zee_K$ over $H$, but an appropriately scaled {\it square\/}
and {\it cube\/} of such $x$-coordinates, respectively.
In all cases, the analogue of the Kronecker-Weber theorem for $K$ is the
following {\it second main theorem of complex multiplication}.
\proclaim
{6.15. Theorem}
The ray class field $H_m$ of conductor $m\Zee_K$ of $K$ is generated
over the Hilbert class field $H$ of $K$ by the values of the Weber function
$f_K$ at the non-zero $m$-torsion points of $\CC/\Zee_K$.
\endproclaim
\noindent
In the non-special cases,
the values of $f_K$ at the $m$-torsion points of $\CC/\Zee_K$
are the $x$-coordinates of the non-zero $m$-torsion points
of the elliptic curve $E_K$ in 6.13. 
For $K=\Que(i)$ and $\Que(\zeta_3)$, one uses squares and cubes of
these coordinates.
In all cases, generating $H_m$ over~$H$ essentially amounts to
computing {\it division polynomials\/} $T_m\in \CC[X]$ 
that have these $x$-coordinates as their roots.
We will define these polynomials as elements of $H[x]$,
as the recursion formulas at the end of this section show that
their coefficients are elements of the ring generated over $\Zee$
by the coefficients of the Weierstrass model of $E_K$.

If $m$ is odd, the non-zero $m$-torsion points come in pairs $\{P, -P\}$
with the same $x$-coordinate $x_P=x_{-P}$,
and we can define a polynomial $T_m(x)\in H[x]$
of degree $(m^2-1)/2$ up to sign by
$$
T_m(x)^2=
m^2 \prod_{P\in E_K[m](\CC)\atop P\ne O} (x-x_P).
$$
For even $m$, we adapt the definition by excluding the
2-torsion points satisfying $P=-P$ from the product,
and define $T_m(X)\in H[x]$ of degree $(m^2-4)/2$ by
$$
T_m(x)^2=
(m/2)^2 \prod_{P\in E_K[m](\CC)\atop 2P\ne O} (x-x_P).
$$
The `missing' $x$-coordinates of the non-zero 2-torsion points are
the zeroes of the cubic polynomial $w_K\in H[x]$ in 6.13.
This is a square in the {\it function field\/}
of the elliptic curve 6.13, and in many ways the natural object to consider is
the `division polynomial'
$$
\psi_m(x,y)=
\cases
T_m(x)      &\text{if $m$ is odd;}\cr
2yT_m(x)     &\text{if $m$ is even.}\cr
\endcases
$$
This is an element of the function field $\CC(x,y)$ living in the
quadratic extension $H[x, y]$ of the polynomial ring $H[x]$
defined by $y^2=w_K(x)$.
It is uniquely defined up to the sign choice we have for $T_m$.
Most modern texts take the sign of the highest coefficient of $T_m$
equal to 1.
Weber [\WE, p.~197] takes it equal to $(-1)^{m-1}$, which amounts
to a sign change $y\mapsto -y$ in $\psi_m(x,y)$.
 
By construction, the function $\psi_m$ has divisor
$(1-m^2)[O]+ \sum_{P\ne O, mP=O} [P]$.
The normalizing highest coefficients $m$ and $m/2$ in $T_m$ lead to
neat recursive formulas
$$
\eqalign{
\psi_{2m+1}&=\psi_{m+2}\psi_m^3-\psi_{m+1}^3\psi_{m-1},\cr
\psi_{2m}  &=(2y)^{-1} \psi_m(\psi_{m+2}\psi_{m-1}^2-\psi_{m+1}^2\psi_{m-2})\cr
}
$$
for $\psi_m$ that are
valid for $m>1$ and $m>2$.
These can be used to compute $\psi_m$ and $T_m$ recursively,
using `repeated doubling' of $m$.
One needs the initial values $T_1=T_2=1$ and
$$
\eqalign{
T_3&=3X^4+ 6aX^2+12bX-a^2,\cr
T_4&=2X^6+10aX^4+40bX^3-10a^2X^2-8abX-16b^2-2a^3,\cr
}
$$
where we have written the Weierstrass polynomial in 6.13 as 
$w_K=4(x^3+ax+b)$ to indicate the relation with the
nowadays more common affine model $y^2=x^3+ax+b$ to which
$E_K$ is isomorphic under $(x,y)\mapsto(x, y/2)$.

\head
7. Class fields from modular functions
\endhead
\noindent
The algorithms in the previous section are based on the main theorems 6.9
and 6.15 of complex multiplication. 
These can be found already
in Weber's 1908 textbook [\WE], so they predate the
class field theory for general number fields.
The oldest proofs of 6.9 and 6.15 are of an analytic nature, 
and derive arithmetic information from
congruence properties of Fourier expansions such as 6.10.
Assuming general class field theory, one can shorten these
proofs as it suffices, just as after 2.2, to show
that, up to sets of primes of zero density,
the `right' primes split completely in the purported class fields.
In particular, it is always possible to restrict attention to the primes
of $K$ of residue degree one in such arguments.
Deuring provides analytic proofs of both kinds in his
survey monograph [\DE].

Later proofs [\LAEF, Part 2] of Deuring and Shimura combine
class field theory with
the {\it reduction\/} of the endomorphisms in 6.6 modulo primes,
which yields endomorphisms of elliptic curves over finite fields.
These proofs are firmly rooted in the algebraic theory of elliptic curves [\SI].
Here one takes for $E_\Lambda$ in 6.6 an elliptic curve $E$ 
that has CM by $\Zee_K$ and is
given by a Weierstrass equation over the splitting field $H'$ over $K$
of the Hilbert class polynomial $\Hil_K(X)$ in 6.8.
In this case the Weierstrass equation can be considered modulo any prime
$\gothq$ of $H'$, and for almost all primes, known as the primes
{\it of good reduction\/}, this yields an elliptic curve
$E_\gothq=(E\mod\gothq)$ over the finite field
$k_\gothq=\Zee_{H'}/\gothq$.
For such $\gothq$, the choice of
an extension of $\gothq$ to $\overline \Que$ yields
a reduction homomorphism $E_K(\overline\Que)\to E_\gothq(\overline k_\gothq)$
on points that is {\it injective\/} on torsion points of order
coprime to $\gothq$.
The endomorphisms of $E$ are given by rational functions 
with coefficients in $H'$, and for primes $\gothq$ of $H'$ of good reduction
there is a natural reduction homomorphism
$$
\End(E)\to\End(E_\gothq)
$$
that is injective and preserves degrees.
The `complex multiplication' by an element
$\alpha\in \End(E)=\Zee_K$ multiplies
the invariant differential $dx/y$ on $E$ by $\alpha$, so it
becomes inseparable in $\End(E_\gothq)$ if and only if
$\gothq$ divides $\alpha$.
The first main theorem of complex multiplication 6.9, which
states that $H'$ equals the Hilbert class field $H$ of $K$
and provides the Galois action on the roots of $\Hil_K(X)$,
can now be derived as follows.
\medskip\noindent
{\bf Proof of 6.9.}
Let $\gothp$ be a prime of degree one of $K$ that is coprime
to the discriminant of $\Hil_K(X)$, and
$\alpha\in\Zee_K$ an element of order 1 at $\gothp$, say
$\alpha\gothp^{-1}=\gothb$ with $(\gothb, \gothp)=1$.
Let $\gotha$ be a fractional $\Zee_K$-ideal.
Then the complex multiplication $\alpha: \CC/\gotha\to\CC/\gotha$
factors in terms of complex tori as
$$
\CC/\gotha \mapright^{\text{can}} \CC/\gotha\gothp^{-1}
\ \smash{\mathop{\longrightarrow}\limits_{\alpha}^\sim}\ 
\CC/\gotha\gothb \mapright^{\text{can}} \CC/\gotha.
$$
If $E$ and $E'$ denote Weierstrass models over $H'$ for
$\CC/\gotha$ and $\CC/\gotha\gothp^{-1}$, we obtain isogenies
$E\to E'\to E$ of degree $p=N\gothp$ and $N\gothb$
with composition $\alpha$.
If we assume that $E$ and $E'$ have good reduction above $\gothp$, 
we can reduce the isogeny $E\to E'$ at some prime $\gothq|\gothp$
to obtain an isogeny $E_\gothq\to E'_\gothq$ of degree $p$.
This isogeny is inseparable as $\alpha$ lies in $\gothq$, and therefore
equal to the Frobenius morphism $E_\gothq\to E_\gothq^{(p)}$
followed by an isomorphism $E_\gothq^{(p)}\isar E'_\gothq$.
The resulting equality $j(E_\gothq^{(p)})=j(E'_\gothq)$ of $j$-invariants
amounts to $j(\gotha)^p=j(\gotha\gothp^{-1})\mod \gothq$, and this implies that
the Frobenius automorphism $\sigma_\gothq\in\Gal(H'/K)$ acts as
$j(\gotha)^{\sigma_\gothq}=j(\gotha\gothp^{-1})$, independently of
the choice of the extension prime $\gothq|\gothp$.
As the $j$-function is an invariant for the homothety class of a lattice,
we have $j(\gotha)=j(\gotha\gothp^{-1})$ if and only if $\gothp$ is
principal.
It follows that up to finitely many exceptions, the primes $\gothp$
of degree one splitting completely in $K\subset H'$ are the principal
primes, so $H'$ equals the Hilbert class field $H$ of $K$, and the
splitting primes in $K\subset H$ are {\it exactly\/} the principal primes.
Moreover, we have $j(\gotha)^{\sigma_\gothc}=j(\gotha\gothc^{-1})$
for the action of the Artin symbol $\sigma_\gothc$, and
$\Hil_K$ is irreducible over $K$ as its roots
are transitively permuted by $\Gal(H/K)=\Cl_K$.  \hfill$\square$
\medskip\noindent
In a similar way, one can understand the content of the second main theorem
of complex multiplication.
If $E_K$ is a Weierstrass model for $\CC/\Zee_K$ defined over the Hilbert
class field $H$ of $K$, then the torsion points in $E_K(\CC)$ have
algebraic coordinates. 
As the group law is given by algebraic
formulas over $H$, the absolute Galois group $G_H$ of $H$ acts
by group automorphisms on $E_K^\tor(\CC)\iso K/\Zee_K$.
Moreover, the action of $G_H$ commutes with the complex multiplication action
of $\End(E_K)\iso\Zee_K$, which is given by isogenies defined over $H$.
It follows that $G_H$ acts by $\Zee_K$-module automorphisms on $E_K^\tor(\CC)$.
For the cyclic $\Zee_K$-module $E_K[m](\CC)\iso {1\over m}\Zee_K/\Zee_K$
of $m$-torsion points, the resulting Galois representation
$$
\textstyle
G_H\tto\Aut_{\Zee_K}({1\over m}\Zee_K/\Zee_K)\iso (\Zee_K/m\Zee_K)^*
$$
of $G_H$ is therefore {\it abelian\/}.
It shows that, just as in the cyclotomic case 6.2,
the Galois action {\it over\/} $H$ on the {\it $m$-division field\/}
of $E_K$, which is the extension of $H$ generated by the $m$-torsion
points of $E_K$, comes from {\it multiplications\/}
on ${1\over m}\Zee_K/\Zee_K$
by integers $\alpha\in\Zee_K$ coprime to $m$.
The content of 6.15 is that, in line with 2.8, we obtain the
$m$-th ray class field of $K$ from this $m$-division field 
by taking invariants under the action of $\Zee_K^*=\Aut(E_K)$.
In the `generic case' where $\Zee_K^*=\{\pm1\}$ has order 2,
adjoining $m$-torsion points `up to inversion' amounts to the equality
$$
H_m=H\big(\{x_P: P\in E_K[m](\CC), P\ne O\}\big)
\leqno(7.1)
$$
occurring in 6.15, since the $x$-coordinate $x_P$ determines $P$ up to
multiplication by $\pm1$.
More generally, a root of unity $\zeta\in\Zee_K^*$ acts
as an automorphism of $E_K$ by $x_{[\zeta]P}=\zeta^{-2} x_P$, so
in the special cases where $K$ equals $\Que(i)$
or $\Que(\zeta_3)$ and $\Zee_K^*$ has order $2k$
with $k=2, 3$, one replaces $x_P$ by $x_P^k$ in 7.1.
The classical Weber functions replacing 6.14 for $K=\Que(i)$ and
$K=\Que(\zeta_3)$ are
$f_K(z)= (g_2^2(\Zee_K)/\Delta(\Zee_K))\wp^2_{\Zee_K}(z)$
and
$f_K(z)= (g_3(\Zee_K)/\Delta(\Zee_K))\wp^3_{\Zee_K}(z)$.
\medskip\noindent
{\bf Proof of 6.15.}
As in the case of 6.9, we show that
the primes of degree one of $K$ that split completely in the extension
$K\subset H'_m$ defined by adjoining to $H$ the $m$-torsion points of $E_K$
`up to automorphisms' are, up to a zero density subset of primes,
the primes in the ray group $R_m$.
Primes $\gothp$ of $K$ splitting in $H'_m$ are principal as they split in $H$.
For each $\Zee_K$-generator $\pi$ of $\gothp$, which is uniquely
determined up to multiplication by $\Zee_K^*$,
one obtains a complex multiplication by $\pi\in \Zee_K\iso\End(E_K)$ that
fixes $E_K[m](\CC)$ `up to automorphisms' if and only if we have 
$\gothp\in R_m$.

Let $\gothp=\pi\Zee_K$ be a prime of degree 1 over $p$ for which
$E_K$ has good reduction modulo $\gothp$.
Then the isogeny $\phi_\pi: E_K\to E_K$, which corresponds to
multiplication by $\pi$ as in 6.5,
reduces modulo a prime $\gothq|\gothp$ of the $m$-division field
of $E_K$ to an endomorphism of degree $p$.
As $\pi$ is in $\gothq$, this reduction is inseparable, so it equals the
Frobenius endomorphism of $E_{K,\gothq}$ up to an automorphism.
One shows [\LAEF, p.\ 125] that this {\it local\/} automorphism
of $E_{K,\gothq}$ is induced by a global automorphism of $E_K$,
i.e., a complex multiplication by a unit in $\Zee_K^*$,
and concludes that $\phi_\pi$ induces a Frobenius automorphism
above $\gothp$ on $H'_m$.
As the reduction modulo $\gothq$ induces an isomorphism
$E_K[m]\isar E_{K,\gothq}[m]$ on the $m$-torsion points,
this Frobenius automorphism is trivial if and only if
we have $\pi\congr1\mmod m\Zee_K$ for a suitable choice of $\pi$.
Thus, $\gothp$ splits completely in $H'_m$ if and only if
$\gothp$ is in $R_m$, and $H'_m$ is the ray class field $H_m$.
\hfill$\square$
\medskip\noindent
The argument just given shows that we have
a concrete realization of the Artin isomorphism
$
(\Zee_K/m\Zee_K)^*/\im[\Zee_K^*] \isar
\Gal(H_m/H) 
$
from 2.8 by complex multiplications.
Passing to the projective limit, this yields the analogue
$$
\widehat\Zee_K^*/\Zee_K^* \isar \Gal(K_\ab/H)\subset G_K^\ab
$$
of 3.1.
For the analogue of 6.3, we note first that for imaginary quadratic
$K$, the subgroup $U_\infty$ in 3.3, which equals $\CC^*$,
maps isomorphically to the connected component of the unit element
in $\AA_K^*/K^*$. As it is the kernel of the Artin map $\psi_K$ in 3.7,
we obtain a commutative diagram
$$
\matrix
\widehat \Zee_K^* &\mapright^{-1}& \widehat \Zee_K^*/\Zee_K^* 
          &\subset &\AA_K^*/(K^*\cdot\CC^*)\cr
\mapdown^{\text{can}}_\wr& &\mapdown_\wr& &\mapdown^{\text{Artin}}_\wr\cr
\Aut_{\Zee_K}(K/\Zee_K)&\tto& \Gal(K_\ab/H)&\subset&\Gal(K_\ab/K)\cr
\endmatrix
\leqno(7.2)
$$
in which the inversion map ``$-1$'' arises just as in 3.9.
A slight difference with the diagram 6.3 for $\Que$ is
that the horizontal arrows now have a small finite kernel coming from
the unit group $\Zee_K^*$. 
Moreover, we have only accounted for the automorphisms of $K_\ab$
over $H$, not over $K$.
Automorphisms of $H_m$ that are not the identity on $H$ arise as
Artin maps $\sigma_\gothc$ of non-principal ideals $\gothc$ coprime to $m$, 
and the proof of 6.9 shows that for the isogeny $\phi_\gothc$ in
the commutative diagram
$$
\matrix
\CC/\Zee_K &\mapright^{\text{can}} & \CC/\gothc^{-1}\cr
\mapdown^W_\wr& &\mapdown^{W'}_\wr\cr
E_K&\mapright^{\phi_\gothc}&E'_K,\cr
\endmatrix
\leqno(7.3)
$$
we have to compute the restriction $\phi_\gothc: E_K[m]\to E'_K[m]$ 
to $m$-torsion points.
To do so in an efficient way, we view
the $j$-values and $x$-coordinates of torsion points involved
as weight zero functions on complex lattices such as $\Zee_K$
or $\gothc$.
As we may scale all lattices as we did for 6.10
to $\Zee\tau+\Zee$ with $\tau\in\HH$,
such functions are {\it modular functions\/} $\HH\to\CC$
as defined in [\LAEF, Chapter 6].

The $j$-function itself is the primordial modular function:
a holomorphic function on $\HH$
that is invariant under the full modular group $\Sl_2(\Zee)$.
Every meromorphic function on $\HH$ that is invariant
under $\Sl_2(\Zee)$ and, when viewed as a function of
$q=e^{2\pi i\tau}$, meromorphic in $q=0$,
is in fact a rational function of $j$.
The Weber function $f_K$ in 6.14 is a function
$$
f_\tau(z)={g_2(\tau)g_3(\tau)\over\Delta(\tau)}\wp_{[\tau,1]}(z)
$$
that depends on the lattice $\Zee_K=\Zee\tau+\Zee=[\tau, 1]$, and 
fixing some {\it choice\/} of a generator $\tau$ of $\Zee_K$ over $\Zee$,
we can label its $m$-torsion values used in generating $H_m$ as
$$
F_u(\tau)=f_\tau(u_1\tau+u_2)\quad
    \text{with }u=(u_1,u_2)\in {1\over m}\Zee^2/\Zee^2\setminus \{(0,0)\}.
\leqno{(7.4)}
$$
For $m>1$, the functions $F_u:\HH\to\CC$ in 7.4 are known as the
{\it Fricke functions\/} of level $m$.
These are holomorphic functions on $\HH$
that are $x$-coordinates of $m$-torsion points on a
`generic elliptic curve' over $\Que(j)$ with $j$-invariant $j$. 
As they are zeroes of division polynomials in $\Que(j)[X]$,
they are algebraic over $\Que(j)$ and
generate a finite algebraic extension of $\Que(j)$,
the $m$-th modular function field
$$
{\Cal F}_m=\Que(j, \{F_u\}_{u\in ({1\over m}\Zee/\Zee)^2\setminus \{(0,0)\}}).
\leqno{(7.5)}
$$
Note that we have ${\Cal F}_1=\Que(j)$ in 7.5.
We may now rephrase the main theorems of complex multiplication
in the following way.
\proclaim
{7.6. Theorem}
Let $K$ be an imaginary quadratic field with ring of integers $\Zee[\tau]$.
Then the $m$-th ray class field extension $K\subset H_m$ is generated
by the finite values $f(\tau)$ of the functions $f\in {\Cal F}_m$.
\endproclaim
\noindent
For generic $K$ this is a directly clear from 6.9 and 6.15,
in the special cases $K=\Que(i), \Que(\zeta_3)$
the functions $F_u(\tau)$ in 7.4 vanish at the generator
$\tau$ of $\Zee_K$, so an extra argument [\LAEF, p. 128] involving
modified Weber functions in ${\Cal F}_m$ is needed.

It is not really necessary to take $\tau$ in 7.6 to be a generator of
$\Zee_K$; it suffices that the elliptic curve $\CC/[\tau, 1]$ is
an elliptic curve having CM by $\Zee_K$.

For computations, it is essential to have the explicit Galois action
of $\Gal(K_\ab/K)$ on the values $f(\tau)$ from 7.6 for
arbitrary functions $f$ in the modular function field
$
{\Cal F}=\cup_{m\ge1} {\Cal F}_m.
$ 
As class field theory gives us the group $\Gal(K_\ab/K)$ in 7.2
as an explicit quotient of the idele class group
$\AA_K^*/K^*$ under the Artin map 3.7,
this means that we need to find the natural action
of $x\in\AA_K^*$ on the values $f(\tau)$ in 7.6.
We will do so by reinterpreting the action of the 
Artin symbol $\sigma_x\in\Gal(K_\ab/K)$ on the function
value of $f$ at $\tau$ as the value of some {\it other\/}
modular function $f^{g_\tau(x)}$ at $\tau$, i.e.,
$$
(f(\tau))^{\sigma_x}= f^{g_\tau(x)}(\tau),
\leqno{(7.7)}
$$
for some natural homomorphism
$
g_\tau: \AA_K^*\to\Aut({\Cal F})
$
induced by $\tau$.

To understand the automorphisms of ${\Cal F}$, we note first that
the natural left action of $\Sl_2(\Zee)$ on $\HH$ gives rise to
a right action on ${\Cal F}_m$ that is easily made explicit 
for the Fricke functions 7.4, using
the `weight 0' property of $f_K$. 
For $M={\alpha\ \beta\atopwithdelims() \gamma\ \delta}\in \Sl_2(\Zee)$ we have
$$
\eqalign{
F_u(M\tau)=F_u\Bigl({\alpha\tau+\beta\over \gamma\tau+\delta}\Bigr)&=
{g_2(\tau)g_3(\tau)\over\Delta(\tau)}
\wp_{[\tau,1]}(u_1(\alpha\tau+\beta)+u_2(\gamma\tau+ \delta))\cr
&=F_{uM}(\tau).\cr }
$$
As $u=(u_1, u_2)$ is in ${1\over m}\Zee^2/\Zee^2$,
we only need to know $M$ modulo $m$, so
the Fricke functions of level $m$ are invariant under
the congruence subgroup
$$
\Gamma(m)=\ker[\Sl_2(\Zee)\to\Sl_2(\Zee/m\Zee)]
$$
of $\Sl_2(\Zee)$, and they are permuted by $\Sl_2(\Zee)$.
As we have $F_{-u_1,-u_2}=F_{u_1, u_2}$, we obtain a natural
right action of $\Sl_2(\Zee/m\Zee)/\{\pm1\}$ on ${\Cal F}_m$.

In addition to this `geometric action', there is a cyclotomic
action of $(\Zee/m\Zee)^*$ on the functions $f\in{\Cal F}_m$ via their
Fourier expansions, which lie in $\Que(\zeta_m)((q^{1/m}))$
as they involve rational expansions in
$$
e^{2\pi i(a_1\tau+a_2)/m} = \zeta_m^{a_2} q^{a_1\over m}
\qquad\text{for $a_1, a_2\in\Zee$.}
$$
On the Fricke function $F_u=F_{(u_1, u_2)}$,
the automorphism
$\sigma_k:\zeta_m\mapsto \zeta_m^k$
clearly induces $\sigma_k: F_{(u_1, u_2)}\mapsto F_{(u_1, ku_2)}$.
Thus, the two actions may be combined to obtain an action of
$\GL_2(\Zee/m\Zee)/\{\pm1\}$ on ${\Cal F}_m$,
with $\Sl_2(\Zee/m\Zee)/\{\pm1\}$ acting geometrically
and $(\Zee/m\Zee)^*$ acting as the subgroup 
$\{\pm{1\ 0\atopwithdelims()0\ k}: k\in (\Zee/m\Zee)^*\}/\{\pm1\}$.
The invariant functions are $\Sl_2(\Zee)$-invariant with
rational $q$-expansion, so they lie in $\Que(j)$ and we have
a natural isomorphism
$$
\GL_2(\Zee/m\Zee)/\{\pm1\}
\isar
\Gal({\Cal F}_m/\Que(j))
\leqno{(7.8)}
$$
or, if we take the union ${\Cal F}=\cup_m{\Cal F}_m$
on the left hand side and the 
corresponding projective limit on the right hand side,
$$
\GL_2(\widehat\Zee)/\{\pm1\}
\isar
\Gal({\Cal F}/\Que(j)).
\leqno{(7.9)}
$$
Note that ${\Cal F}_m$ contains $\Que(\zeta_m)(j)$
as the invariant field of $\Sl_2(\Zee/m\Zee)/\{\pm 1\}$,
and that the action on the subextension $\Que(j)\subset \Que_\ab(j)$
with group $\Gal(\Que_\ab/\Que)\iso\widehat\Zee^*$ is via
the determinant map $\det: \GL_2(\widehat\Zee)/\{\pm1\}\to \widehat\Zee^*$.

In order to discover the explicit form of the homomorphism
$g_\tau$ in 7.7, 
let $\gothp=\pi\Zee_K$ be a principal prime of $K$.
Then the Artin symbol $\sigma_\gothp$ is the identity on $H$,
and the proof of 6.15 shows that its action on the $x$-coordinates
of the $m$-torsion points of $E_K$ for $m$ not divisible
by $\pi$ can be written as
$$
F_u(\tau)^{\sigma_\gothp}
=F_u(\pi\tau)=F_{uM_\pi}(\tau),
$$
where $M_\pi$ is the matrix in $\GL_2(\Zee/m\Zee)$ that represents
the multiplication by $\pi$ on ${1\over m}\Zee_K/\Zee_K$ with respect
to the basis $\{\tau, 1\}$.
In explicit coordinates, this means that if 
$\tau\in\HH$ is a zero of the polynomial $X^2+BX+C$ of discriminant
$B^2-4C=\Delta_K$ and
$\pi=x_1\tau+x_2$ is the representation
of $\pi$ on the $\Zee$-basis $[\tau, 1]$ of $\Zee_K$,
then we have
$$
M_\pi=
\left(
\matrix
-Bx_1+x_2&-Cx_1\cr
 x_1      & x_2\cr
\endmatrix
\right)\in \GL_2(\Zee/m\Zee).
\leqno{(7.10)}
$$
As the Fricke functions of level $m$ generate ${\Cal F}_m$,
we obtain in view of 7.8 the identity
$$
f(\tau)^{\sigma_\gothp}= f^{M_\pi}(\tau)
         \qquad\text{for $f\in{\Cal F}_m$ and $\pi\nmid m$,}
$$
which is indeed of the form 7.7.
We can rewrite this in the style of the diagram 7.2
by observing that the Artin symbol of 
$\pi\in K_\gothp^*\subset \AA_K^*$
acts as $\sigma_\gothp$ on torsion points of order $m$
coprime to $\gothp$, and trivially on $\pi$-power torsion points.
Moreover, $(\pi\mod K^*)\in \AA_K/K^*$
is in the class of the idele $x\in\widehat\Zee_K^*$
having component $1$ at $\gothp$ and $\pi^{-1}$ elsewhere.
Thus, if we define 
$$
\eqalign{
g_\tau: \widehat\Zee_K^*&\tto \GL_2(\widehat\Zee)\cr
                       x&\longmapsto M_x^{-1}\cr}
\leqno{(7.11)}
$$
by sending $x=x_1\tau+x_2\in \widehat\Zee_K$ to the {\it inverse\/}
of the matrix $M_x$
describing multiplication by $x$ on $\widehat\Zee_K$ with respect to the
basis $[\tau, 1]$, then $M_x$ is given explicitly as in 7.10,
and formula 7.7 holds for $f\in{\Cal F}$ and $x\in \widehat\Zee_K^*$
if we use the natural action of $\GL_2(\widehat\Zee)$ on
${\Cal F}$ from 7.9.

To obtain complex multiplication by arbitrary ideles, we
note that on the one hand,
the idele class quotient $\AA_K^*/(K^*\cdot\CC^*)$ from 7.2,
which is isomorphic to $\Gal(K_\ab/K)$ under the Artin map, 
is the quotient of the unit group $\widehat K^*$ of the 
{\it finite adele ring\/}
$$
\widehat K = 
             \widehat\Zee_K\otimes_\Zee \Que =
             \prod\nolimits'_{\gothp\text{ finite}} K_\gothp
\subset \AA_K=\widehat K\times\CC 
$$
by the subgroup $K^*\subset \widehat K^*$ of principal ideles.
On the other hand, not all automorphisms of $\Cal F$ come from 
$\GL_2(\widehat\Zee)$ as in 7.9:
there is also an action of the projective linear group
$\PGL_2(\Que)^+=\GL_2(\Que)^+/\Que^*$ of rational matrices of
positive determinant,
which naturally act on $\HH$ by linear fractional transformations.
It does not fix $j$, as it maps the elliptic curve $\CC/[\tau, 1]$
defined by $\tau\in\HH$ not to an isomorphic, but to an isogenous 
curve.
More precisely, if we pick 
$M={\alpha\ \beta\atopwithdelims() \gamma\ \delta}\in \GL_2(\Que)^+$
in its residue class modulo $\Que^*$ such that
$M^{-1}$ has integral coefficients, then
the lattice 
$$
\textstyle
(\gamma\tau+\delta)^{-1}[\tau, 1]=
[{\tau\over \gamma\tau+\delta}, {1\over \gamma\tau+\delta}]=
M^{-1}[{\alpha\tau+\beta\over \gamma\tau+\delta}, 1]
$$
is a sublattice of finite index $\det M^{-1}$ in 
$
[{\alpha\tau+\beta\over \gamma\tau+\delta}, 1],
$
and putting $\mu=(\gamma\tau+\delta)^{-1}$, we have a commutative diagram
$$
\matrix
\CC/[\tau,1]&\mapright^\mu &\CC/[M\tau,1]
                   \rlap{$=\CC/[{\alpha\tau+\beta\over\gamma\tau+\delta},1]$}\cr
\mapdown^W_\wr& &\mapdown^{W'}_\wr\cr
E_\tau&\mapright^{\phi_\mu}&E_{M\tau}\cr
\endmatrix
$$
as in 7.3. Moreover, the torsion point $u_1\tau + u_2$ having coordinates
$u=(u_1, u_2)$ with respect to $[\tau,1]$
is mapped to the torsion point with coordinates $uM^{-1}$
with respect to the basis $[M\tau, 1]$.

We let $\widehat\Que=\widehat\Zee\otimes_\Zee \Que=\prod'_p \Que_p$
be the ring of finite $\Que$-ideles.
Then every element in the ring $\GL_2(\widehat\Que)$ can be written
as $UM$ with $U\in\GL_2(\widehat\Zee)$ and $M\in \GL_2(\Que)^+$.
This representation is not unique as $\GL_2(\widehat\Zee)$ and
$\GL_2(\Que)^+$ have non-trivial intersection $\Sl_2(\Zee)$,
but we obtain a well-defined action
$\GL_2(\widehat\Que)\to \Aut({\Cal F})$ by putting
$$
f^{UM}(\tau)=f^U(M\tau).
$$
We now extend, for the zero $\tau\in\HH$ of a polynomial $X^2+BX+C\in\Que[X]$,
the map $g_\tau$ in 7.11 to
$$
\eqalign{
g_\tau: \widehat K^*=(\widehat\Que\tau+\widehat\Que)^*
                                  &\tto \GL_2(\widehat\Que)\cr
        x=x_1\tau+x_2\quad &\longmapsto M_x^{-1}=
        \left(\matrix -Bx_1+x_2&-Cx_1\cr x_1   & x_2\cr \endmatrix \right)^{-1}
\cr}
\leqno{(7.12)}
$$
to obtain the complete Galois action of $\Gal(K_\ab/K)\iso \widehat K^*/K^*$
on modular function values $f(\tau)$.
The result obtained is known as {\it Shimura's reciprocity law}.
\proclaim
{7.13. Theorem}
Let $\tau\in\HH$ be imaginary quadratic, $f\in{\Cal F}$ a modular function
that is finite at $\tau$ and $x\in\widehat K^*/K^*$ a finite idele
for $K=\Que(\tau)$.
Then $f(\tau)$ is abelian over $K$, and the idele $x$ acts on it via its
Artin symbol by
$$
f(\tau)^x=f^{g_\tau(x)}(\tau),
$$
where $g_\tau$ is defined as in $7.12$.
\endproclaim

\head
8. Class invariants
\endhead
\noindent
Much work has gone into algorithmic improvements
of the classical algorithms in section 6, with the aim of
reducing the size of the class polynomials obtained.
Clearly the {\it degree\/} of the polynomials involved
cannot be lowered, as these are the degrees of the field extensions
one wants to compute.
There are however methods to reduce the size of their coefficients.
These already go back to Weber [\WE], who made extensive
use of `smaller' functions than $j$ to compute class fields
in his 1908 algebra textbook.
The function $\gothf_2$ that we used to compute $j$ in 6.12, and that
carries Weber's name (as does the elliptic function in 6.14) provides
a good example.
A small field such as $K=\Que(\sqrt{-71})$, for which the class group
of order 7 is easily computed by hand, already has the sizable
Hilbert class polynomial
        $$\eqalign{
        \text{Hil}_K(X) &=  \;X^7 + 313645809715\;X^6 -
         3091990138604570\;X^5 \cr&\quad
        + 98394038810047812049302\;X^4
        - 823534263439730779968091389\;X^3 \cr &\quad
        + 5138800366453976780323726329446\;X^2 \cr&\quad
        - 425319473946139603274605151187659\; X \cr &\quad
        + 737707086760731113357714241006081263 .\cr
        }$$
However, the Weber function $\gothf_2$, when evaluated at an appropriate
generator of $\Zee_K$ over~$\Zee$, also yields a
generator for $H$ over $K$, and its irreducible polynomial is
$$
X^7+X^6-X^5-X^4-X^3+X^2+2X+1.
$$
As Weber showed, the function $\gothf_2$ can be used to generate $H$
over $K$ when 2 splits and 3 does not ramify in $\Que\subset K$.
The general situation illustrated by this example
is that, despite the content of 7.6,
it is sometimes possible to use a function $f$ of
high level, like the Weber function $\gothf_2\in{\Cal F}_{48}$
of level 48, to generate the Hilbert class field $H$ of conductor 1.
The attractive feature of such high level functions $f$ is that they
can be much smaller than the $j$-function itself.
In the case of $\gothf_2$, the extension $\Que(j)\subset\Que(\gothf_2)$
is of degree 72 by 6.12, and this means that the size of the coefficients
of class polynomials using $\gothf_2$ is about a factor 72 smaller than 
the coefficients of $\Hil_K(X)$ itself.
Even though this is only a constant factor, and complex multiplication
is an intrinsically `exponential' method, the computational improvement
is considerable.
For this reason,
Weber's use of `small' functions has gained renewed interest in
present-day computational practice. 

Shimura's reciprocity law 7.13 is a convenient tool to understand 
the occurrence of {\it class invariants\/}, i.e., modular functions
$f\in{\Cal F}$ of higher level that generate the Hilbert class field
of $K$ when evaluated at an appropriate generator $\tau$ of $\Zee_K$.
Classical examples of such functions used by Weber are
$\gamma_2=\root 3\of j$ and $\gamma_3=\sqrt{j-1728}$, which have level 3 and 2.
As is clear from 6.12, the $j$-function can also be constructed out of
even smaller building blocks involving the Dedekind $\eta$-function 6.11.
Functions that are currently employed in actual computations are
$$
{\eta(pz)\over \eta(z)}\qquad\text{and}\qquad
{\eta(pz)\eta(qz)\over \eta(pqz)\eta(z)},
\leqno{(8.1)}
$$
which are of level $24p$ and $24pq$.
These functions, or sometimes small powers of them,
can be used to generate $H$, and the resulting minimal polynomials
have much smaller coefficients than $\Hil_K(X)$.
We refer to [\COB, Section 6.3] for the precise theorems,
and indicate here how to use 7.13 to obtain such results
for arbitrary modular functions $f\in{\Cal F}$.

Let $f\in{\Cal F}$ be any modular function of level $m$, and
assume $\Que(f)\subset {\Cal F}$ is Galois.
Suppose we have an explicit Fourier expansion in $\Que(\zeta_m)((q^{1/m}))$
that we can use to approximate its values numerically.
Suppose also that we know the explicit action of the generators
$S: z\mapsto 1/z$ and $T: z\mapsto z+1$ on $f$.
Then we can determine the Galois orbit of $f(\tau)$ for an element
$\tau\in\HH$ that generates $\Zee_K$ in the following way.
First, we determine elements $x=x_1\tau+x_2\in \Zee_K$
with the property that they generate $(\Zee_K/m\Zee_K)^*/\Zee_K^*$.
Then the Galois orbit of $f(\tau)$ over $H$ is determined using 7.13,
and amounts to computing the (repeated) action of the
matrices $g_\tau(x)\in\GL_2(\Zee/m\Zee)$
(given by the right hand side of 7.10) on $f$.
This involves writing $g_\tau(x)$ as a product of powers of $S$ and $T$
and a matrix ${1\ 0\atopwithdelims()0\ k}$ acting on $f$ via its Fourier
coefficients.
Although $f$ may have a large $\GL_2(\Zee/m\Zee)$-orbit over $\Que(j)$,
the matrices $g_\tau(x)$ only generate a small subgroup of $\GL_2(\Zee/m\Zee)$
isomorphic to $(\Zee_K/m\Zee_K)^*/\Zee_K^*$, and one often finds that the orbit
of $f$ under this subgroup is quite small.
In many cases, one can slightly modify $f$, multiplying it by suitable
roots of unity or raising it to small powers, to obtain an orbit
of length one.
This means that $f\in{\Cal F}$ is invariant under
$g_\tau[\widehat\Zee_K^*]\subset\GL_2(\widehat\Zee)$.
As we have the fundamental equivalence
$$
f(\tau)^x=f(\tau)\quad\Leftrightarrow\quad
f^{g_\tau(x)}=f,
\leqno{(8.2)}
$$
this is equivalent to finding that $f(\tau)$ is a class invariant
for $K=\Que(\tau)$.
The verification that $g_\tau[\widehat\Zee_K^*]$ stabilizes $f$
takes place modulo the level $m$ of $f$, so it follows from 7.12
that if $f(\tau)$ is a class invariant
for $K=\Que(\tau)$, then $f(\tau')$ is a class invariant
for $K'=\Que(\tau')$ whenever $\tau'\in\HH$ is a generator of $\Zee_{K'}$
that has an irreducible polynomial congruent modulo $m$
to that of $\tau$.
In particular, a function of level $m$ that yields class invariants does
so for families of quadratic fields for which the discriminant is in
certain congruence classes modulo~$4m$.

If $f(\tau)$ is found to be a class invariant, we need to
determine its conjugates over $K$ to determine its
irreducible polynomial over $K$ we did in 6.8 for $j(\tau)$.
This amounts to computing $f(\tau)^{\sigma_\gothc}$ as in 6.9,
with $\gothc$ ranging over the ideal classes of $\Cl_K$.
If we list the ideal classes of $\Cl_K$ as in Section 6
as integer triples $(a, b, c)$ representing the reduced quadratic forms
of discriminant $\Delta_K$, the Galois action of their Artin symbols 
in 6.9 may be given by
$$
\textstyle
j(\tau)^{(a, -b, c)}= j({-b+\sqrt{b^2-4ac}\over 2a}).
$$
For a class invariant $f(\tau)$ a similar formula is provided
by Shimura's reciprocity law. 
Let $\gotha=\Zee\cdot {-b+\sqrt{b^2-4ac}\over 2}+\Zee\cdot a$
be a $\Zee_K$-ideal in the ideal class corresponding to the
form $(a, b, c)$.
Then the $\widehat\Zee_K$-ideal $\gotha\widehat\Zee_K$ is
principal as $\Zee_K$-ideals are locally principal, and we
let $x\in \widehat \Zee_K$ be a generator.
The element $x$ is a finite idele in $\widehat K^*$, and the Artin
symbol of $x^{-1}$ acts on $f(\tau)$ as the Artin symbol
of the form $(a, -b, c)$.
We have $U=g_\tau(x^{-1}) M^{-1}\in \GL_2(\widehat\Zee)$ for the
matrix $M\in\GL_2(\Que)^+$ defined by
$[\tau, 1]M=[{b+\sqrt{b^2-4ac}\over 2}, 2a]$, since $U$ stabilizes
the $\widehat\Zee_K$-lattice spanned by the basis $[\tau, 1]$.
Applying 7.13 for the idele $x^{-1}$ yields the desired formula
$$
\textstyle
f(\tau)^{(a, -b, c)}= f^U({-b+\sqrt{b^2-4ac}\over 2a}).
$$
This somewhat abstract description may be phrased as a
simple explicit recipe for the coefficients of $U\in \GL_2(\widehat\Zee)$,
which we only need to know modulo $m$, see [\STB].

There are limits to the improvements coming from intelligent choices of
modular functions to generate class fields.
For any non-constant function $f\in{\Cal F}$, there is a polynomial
relation $\Psi(j, f)=0$ between $j$ and $f$, with
$\Psi\in\CC[X, Y]$ some irreducible polynomial with algebraic coefficients.
The {\it reduction\/} factor one obtains by using class invariants coming from 
$f$, if these exist, instead of the classical $j$-values is defined as
$$
r(f)=
{\deg_f(\Psi(f,j))\over \deg_j(\Psi(f,j))}.
$$
By [\HS, Proposition B.3.5], this is, asymptotically, the {\it inverse\/}
of the factor
$$
\lim_{h(j(\tau)) \rightarrow \infty}
{h(f(\tau)) \over h(j(\tau))}.
$$
Here $h$ is the absolute logarithmic height, and we take the limit over
all CM-points $\Sl_2(\Zee)\cdot\tau\in\HH$.
It follows from gonality estimates for modular curves
[\BS, Theorem 4.1] that $r(f)$ is bounded from above by
$24/\lambda_1$, where $\lambda_1$ is a lower bound for the smallest
positive eigenvalue of the Laplacian (as defined in [\SA])
on modular curves.
The best proved lower bound $\lambda_1\ge {975\over 4096}$
[\KI, p. 176] yields $r(f)\le 32768/325\approx 100.8$,
and Selberg's conjectural bound $\lambda_1\ge 1/4$ implies $r(f)\le 96$.
Thus Weber's function $\gothf_2$, which has $r(f)=72$ and
yields class invariants for a positive density subset of
all discriminants, is close to being optimal.
\medskip\noindent
{\bf Acknowledgements.}
Useful comments on earlier versions of this paper were provided
by Reinier Br\"oker, Ren\'e Schoof and Marco Streng.
Bjorn Poonen provided us with the reference to [\KI].

\enddocument